\documentclass[12pt]{article}
\usepackage{graphicx}
\usepackage{amssymb}
\usepackage{amsmath}
\numberwithin{equation}{section}

\begin{document}
\author{Lev Sakhnovich}
\date{May 15, 2006}
\textbf{Integral Equations in the Theory of Levy Processes}

\begin{center} Sakhnovich Lev \end{center}
735 Crawford ave., Brooklyn, 11223, New York, USA \\
 E-mail
address: lev.sakhnovich@verizon.net
\begin{center}Abstract \end{center}
In this article we consider the Levy processes and the corresponding
semigroup. We represent the generator of this semigroup  in a
convolution form. Using the obtained convolution form and the theory
of integral equations we investigate the properties of a wide class
of Levy processes  (potential, quasi-potential, the probability of
the Levy process remaining within the given domain, long time
behavior, stable processes). We analyze in detail a number of
concrete examples of the Levy processes (stable  processes, the
variance damped Levy processes, the variance gamma processes, the
normal Gaussian process, the Meixner process, the compound Poisson
process).\\
 \textbf{Mathematic Subject Classification (2000).}
Primary 60G51,
Secondary 60J45, 60G17, 45A05.

\begin{center}{Introduction}\end{center}
In the famous article by M.Kac [11] a number of examples demonstrate
the interconnection between the probability theory and the theory of
integral and differential equations. In particular, in article [11] the
Cauchy process was considered. Later M.Kac method was used both for
symmetric stable processes [30], [19] and non-symmetric stable
processes [20]-[22]. In the present article with the help of
M.Kac's idea [11] and the theory of integral equations with the
difference kernels [22] we investigate a wide class of Levy
processes. We note that within the last ten years the Levy processes
have got a number of new important applications, particularly to
financial problems. We consider separately the examples of Levy
processes which are used in financial mathematics.\\
Now we  shall formulate the main results of the article.\\
1. The Levy process $X_{t}$ defines a strongly continuous semigroup
$P_{t}$ (see [23]). The generator $L$ of the semigroup $P_{t}$ is a
pseudo-differential operator. We show that for a broad class of the
Levy processes  the generator $L$ can be represented in a
convolution type form (section 2):
\begin{equation}
Lf=\frac{d}{dx}S\frac{d}{dx}f,\end{equation} where the operator
$S$ is defined by the relation
\begin{equation}
Sf=\frac{1}{2}Af+\int_{-\infty}^{\infty}k(y-x)f(y)dy,\end{equation}
Such a representation opens a possibility to use the theory of integral equations
with difference kernels [22].\\
2. We introduce the notion of the truncated
generator $L_{\Delta}$, which is important in our theory and coincides with $L$ on the bounded
system of segments $\Delta$. We define the quasi-potential
 $B$  by the relation $-L_{\Delta}Bf=f$ (sections 3 and 4). Under our assumptions
 the operator \begin{equation}Bf=\int_{\Delta}\Phi(x,y)f(y)dy \end{equation}  is compact.
  It means that the operator $B$ has a discrete spectrum
 $\lambda_{j}\quad (j=1,2,...),\quad \lambda_{j}{\to}0.$
 Hence, the corresponding truncated generator $L_{\Delta}$ has a discrete spectrum
 too. Using the results from the theory of the integral equations with
  the difference kernels we represent the quasi-potential $B$ in the explicit form.\\
3.The probability $p(t,\Delta)$ of the Levy process remaining within
 the given domain $\Delta$
 (ruin problem) is investigated in our paper
(section 5). M.Kac [11] had obtained the first results of this type
for symmetric stable processes. Later we transposed these results
for the non-symmetric stable processes [20], [22]. In this paper we
show that integro-differential equation of M.Kac type is true for
all the Levy processes, which have a continuous density. M.Kac
writes [11]: "We are led here to integro-differential equations
which offer formidable analytic difficulties and which we were able
to solve only in very few cases." M. Kac was able to overcome these
difficulties only for Cauchy processes. H.Widom [30] solved these
equations for symmetric stable processes. Both symmetric and
non-symmetric stable processes were investigated in our works
[19]-[22].
Now we develop these results and transfer them on the wide class of the Levy processes.\\
4. In sections 6-8 we investigate the structure and the properties
of the  quasi-potential $B$. In particular, we prove that the
operator $B$ is compact and the following important inequality:
\begin{equation}\Phi(x,y){\geq}0\end{equation}
is true.  From inequality (0.4) and Krein-Rutman theorem [13] we
deduce that the operator $B$ has a positive eigenvalue $\lambda_{1}$, which
is no
less in modulus than every other eigenvalue of $B$.\\
 5. Section 9
contains formulas and estimations for the probability $p(t,\Delta)$
that a
 sample of the process $X_{\tau}$ remains inside the given domain
  $\Delta$, when
 $0{\leq}\tau{\leq}t$. Under certain conditions we have  obtained the
 asymptotic formula
\begin{equation}p(t,\Delta)=e^{-t/\lambda_{1}}[c_{1}+o(1)],\quad t{\to}\infty.\end{equation}
 6. In sections 10-12 we separately consider a
special case of the Levy processes, i.e., stable processes. We use
the notation $p(t,a)=p(t,\Delta)$ when $\Delta=[-a,a].$ We consider
the important case, where $a$ depends on $t$ and $a(t){\to}\infty$
when $t{\to}\infty.$ We compare the obtained results with well-known
results (the iterated logarithm law, the results for the first
hitting time, the results for the most visited sites problems).
Further we introduce the notation $p(t,-b,a)=p(t,\Delta)$ when
$\Delta=[-b,a].$ In case of the Wiener process we found the
asymptotic behavior of $p(t,-b,a)$ when $b{\to}\infty.$ It is easy
to see that $p(t,-\infty,a)$ coincides with the formula for the
first hitting time.\\
 7. We analyze in detail a number of
concrete examples of the Levy processes which are used in the
financial mathematics (stable processes, the variance damped Levy
processes, the variance gamma processes, the normal Gaussian
process,
the Meixner process, compound Poisson process.)\\
\section{Main notions}
 Let us consider the Levy processes $X_{t}$ on $R$. If
$P(X_{0}=0)=1$ then Levy-Khinchine formula gives (see[4],[23])
\begin{equation}
\mu(z,t)=E\{\mathrm{exp}[izX_{t}]\}=
\mathrm{exp}[-t\lambda(z)],\quad t{\geq}0, \end{equation} where
\begin{equation}
\lambda(z)=\frac{1}{2}Az^{2}-i{\gamma}z-\int_{-\infty}^{\infty}(e^{ixz}-1-ixz1_{D(x)})d\nu(x).
\end{equation} Here $A{\geq}0,\quad \gamma=\overline{\gamma},$ and
$D=\{x:|x|{\leq}1\}$ is the segment [-1,1], $\nu(x)$ is a
monotonically increasing function  satisfying the conditions
\begin{equation}
\int_{-\infty}^{\infty}\frac{x^{2}}{1+x^{2}}d\nu(x)<\infty.
\end{equation} By $P_{t}(x_{0},\Delta)$ we denote the probability
$P(X_{t}{\in}\Delta)$ when $P(X_{0}=x_{0})=1$ and $\Delta{\in}R$.
The transition operator is defined by the formula \begin{equation}
P_{t}f(x)=\int_{-\infty}^{\infty}P_{t}(x,dy)f(y).\end{equation} Let
$C_{0}$ be the Banach space of continuous functions $f(x) \quad
(-\infty<x<\infty)$ satisfying the condition
$\mathrm{lim}f(x)=0,\quad |x|{\to}\infty$  with the norm $||f||=
\mathrm{sup}_{x}|f(x)|$. We denote by $C_{0}^{n}$ the set of
$f(x){\in}C_{0}$ such that $f^{(k)}(x){\in}C_{0},\quad
(1{\leq}k{\leq}n).$ It is known that [23]
\begin{equation} P_{t}f{\in}C_{0},\end{equation}
if $f(x){\in}C_{0}.$\\
Now we formulate the following important result (see [4],[23]) .\\
\textbf{Theorem 1.1.}\emph{The family of the operators $P_{t}\quad
(t{\geq}0)$ defined by the Levy process $X_{t}$ is a strongly
continuous  semigroup on $C_{0}$ with the norm $||P_{t}||=1$. Let
$L$ be its infinitesimal generator. Then}
\begin{equation} Lf=\frac{1}{2}A\frac {d^{2}f}{dx^{2}}+\gamma
\frac{df}{dx}+\int_{-\infty}^{\infty}(f(x+y)-f(x)-y\frac{df}{dx}1_{D(x)})d\nu(x),\end{equation}
\emph{where} $f{\in}C_{0}^{2}$.
\section{Convolution type form of infinitesimal generator}
1. In this section we prove that under some conditions the
infinitesimal generator $L$ can be represented in the special
convolution type form \begin{equation}
Lf=\frac{d}{dx}S\frac{d}{dx}f,\end{equation} where the operator
$S$ is defined by the relation \begin{equation}
Sf=\frac{1}{2}Af+\int_{-\infty}^{\infty}k(y-x)f(y)dy,\end{equation}
and for arbitrary $M (0<M<\infty)$ we have
\begin{equation} \int_{-M}^{M}|k(t)|dt<\infty.\end{equation}
The representation $L$ in form (2.1) is convenient as the operator
$L$ is expressed with the help of the classic differential and
convolution operators.\\By $C_{c}$ we denote the set of functions
$f(x){\in}C_{0}$ with a compact support.\\ \textbf{Lemma 2.1.}
\emph{Let the following conditions be
fulfilled.\\
1. The function $\nu(x)$ is monotonically increasing, has the
derivative when $x{\ne}0$ and}
\begin{equation}
\int_{-\infty}^{\infty}\frac{x^{2}}{1+x^{2}}d\nu(x)=
\int_{-\infty}^{\infty}\frac{x^{2}}{1+x^{2}}{\nu}^{\prime}(x)dx<\infty,
\end{equation}
\begin{equation}\nu(x){\to}0, \quad x{\to}\infty.\end{equation}

 2.\emph{ For arbitrary $M \quad (0<M<\infty)$ we
have}
\begin{equation}
\int_{-M}^{M}|\nu(x)|dx<\infty,\quad
\int_{-M}^{M}|x|{\nu}^{\prime}(x)dx<\infty.\end{equation} 3.
\begin{equation}
x\nu(x){\to}0,\quad x{\to}0.\end{equation} \emph{Then the expression
\begin{equation}J=\int_{-\infty}^{\infty}[f(y+x)-f(x)]{\nu}^{\prime}(y)dy
\end{equation} can be represented in the convolution type form}
\begin{equation}
J=\frac{d}{dx}\int_{-\infty}^{\infty}f^{\prime}(y)k(y-x)dy
  \end{equation} \emph{ where} $f(x){\in}C_{0}^{2},\quad
k(x)=\int_{0}^{x}\nu(y)dy$.\\ \emph{Proof.} For every
$f(x){\in}C_{c}$  there exists such $M \quad (0<M<\infty)$ that
\begin{equation} f(x)=0 ,\quad x{\notin}[-M,M].\end{equation}Let
us introduce the following notations \begin{equation}
J_{1}=\frac{d}{dx}\int_{-\infty}^{x}f^{\prime}(y)k(y-x)dy,\end{equation}
\begin{equation}
J_{2}=\frac{d}{dx}\int_{x}^{\infty}f^{\prime}(y)k(y-x)dy.\end{equation}
Using (2.11) we have
\begin{equation}J_{1}=
-\frac{d}{dx}\int_{-M}^{x}[f(y)-f(x)+f(x)]k^{\prime}(y-x)dy.\end{equation}
From (2.11) and (2.13) we deduce the relation \begin{equation}
J_{1}=f(x)k^{\prime}(-M-x)+\int_{-M}^{x}[f(y)-f(x)]k^{\prime\prime}(y-x)dy.\end{equation}
When $M{\to}\infty$ we obtain the equality
\begin{equation}
J_{1}=\int_{-\infty}^{0}[f(y+x)-f(x)]k^{\prime\prime}(y)dy.\end{equation}
In the same way we deduce the relation \begin{equation}
J_{2}=\int_{0}^{\infty}[f(y+x)-f(x)]k^{\prime\prime}(y)dy.
\end{equation} Relation (2.9) follows directly from formulas (2.15) , (2.16)
and the equality
 $J=J_{1}+J_{2}.$
 The lemma is proved.\\
\textbf{Lemma 2.2.} \emph{Let the following conditions be
fulfilled.\\1. The function $\nu(x)$ satisfies conditions $(2.4)$
and $(2.5)$ of Lemma
$2.1$.\\
2. For arbitrary $M\quad(0<M<\infty)$ we have
\begin{equation}
\int_{-M}^{M}|k(x)|dx<\infty,\quad
\int_{-M}^{M}|x{\nu}(x)|dx<\infty,\end{equation}where}
\begin{equation}k^{\prime}(x)=\nu(x),\quad x{\ne}0.\end{equation}
 3.
\begin{equation}
xk(x){\to}0,\quad x{\to}0;\quad x^{2}\nu(x){\to}0,\quad
x{\to}0.\end{equation} \emph{Then the equality}
\begin{equation}
J=\int_{-\infty}^{\infty}[f(y+x)-f(x)-y\frac{df(x)}{dx}1_{D(y)}]{\nu}^{\prime}(y)dy
+{\Gamma}f^{\prime}(x),
\end{equation}
\emph{is true, where $\Gamma=\overline{\Gamma}$ and}
$f(x){\in}C_{c}$.\\
\emph{Proof.} From (2.11)  we obtain  the relation
\begin{equation}
J_{1}=f^{\prime}(x)\gamma_{1}-\int_{x-1}^{x}[f^{\prime}(y)-f^{\prime}(x)]k^{\prime}(y-x)dy-
\int_{-M}^{x-1}f^{\prime}(y)k^{\prime}(y-x)dy,\end{equation} where
$\gamma_{1}=k(-1).$ We introduce the notations \begin{equation}
P_{1}(x,y)=f(y)-f(x)-(y-x)f^{\prime}(x),\quad
P_{2}(x,y)=f(y)-f(x).\end{equation}Integrating by parts (2.21) and
passing to the limit when $M{\to}\infty$ we deduce that
\begin{equation}J_{1}=f^{\prime}(x)\gamma_{2}+\int_{x-1}^{x}P_{1}(x,y)k^{\prime\prime}(y-x)dy
+\int_{-M}^{x-1}P_{2}(x,y)k^{\prime\prime}(y-x)dy,\end{equation}
where $\gamma_{2}=k(-1)-k^{\prime}(-1).$ It follows from (2.22) and (2.23)
that \begin{equation}
J_{1}=\int_{-\infty}^{x}[f(y+x)-f(x)-y\frac{df(x)}{dx}1_{D(y)}]{\nu}^{\prime}(y)dy
+{\gamma}_{2}f^{\prime}(x).\end{equation}In the same way it can be
proved that \begin{equation}
J_{2}=\int_{x}^{\infty}[f(y+x)-f(x)-y\frac{df(x)}{dx}1_{D(y)}]{\nu}^{\prime}(y)dy
+{\gamma}_{3}f^{\prime}(x),\end{equation}where
$\gamma_{3}=-k(1)+k^{\prime}(1).$ The relation (2.20) follows directly
from  (2.24) and (2.25). Here $\Gamma=\gamma_{2}+\gamma_{3}.$ The lemma is
proved.\\
\textbf{Remark 2.1.} The operator $L_{0}f=\frac{d}{dx}f$ can be
represented in  form (2.1), (2.2), where \begin{equation}
S_{0}f=\int_{-\infty}^{\infty}p_{0}(y-x)f(y)dy,\end{equation}
\begin{equation}p_{0}(x)=\frac{1}{2}\,\mathrm{sign}(x).\end{equation}\\
From Lemmas 2.1 , 2.2 and Remark 2.1 we deduce the following
assertion.\\
\textbf{Theorem 2.1.} \emph{Let the conditions of Lemma $2.1$ or
Lemma $2.2$ be fulfilled . Then the corresponding operator $L$ has a
convolution type
form }(2.1),(2.2).\\
\textbf{Proposition 2.1.} \emph{The generator $L$ of the Levy
process $X_{t}$ admits the convolution type representation
$(2.1),(2.2)$ if there exist such $C>0$ and $0<\alpha<2,\quad
\alpha{\ne}1$ that}
\begin{equation}\nu^{\prime}(y){\leq}C|y|^{-\alpha-1},\quad
\end{equation} \emph{Proof.} The function $\nu(y)$  has the form
\begin{equation}
\nu(y)=\int_{-\infty}^{y}{\nu}^{\prime}(t)dt1_{y<0}-\int_{y}^{\infty}{\nu}^{\prime}(t)dt1_{y>0}.
\end{equation}First we shall consider the case, when $1<\alpha<2$, and introduce the function
\begin{equation}
k_{0}(y)=\int_{-\infty}^{y}(y-t){\nu}^{\prime}(t)dt1_{y<0}-\int_{y}^{\infty}(y-t){\nu}^{\prime}(t)dt1_{y>0}.
\end{equation}We obtain the relation
\begin{equation}
k(y)=k_{0}(y)+(\gamma-{\Gamma})\mathrm{p}_{0}(y),\quad
1<\alpha<2,\end{equation} where $k_{0}(y)$ and $p_{0}(y)$ are
defined by (2.27) and (2.30) respectively.
 The constant $\Gamma$ is defined by the relation:
\begin{equation}\Gamma=k_{0}(-1)-k_{0}^{\prime}(-1)-k_{0}(-1)+k_{0}^{\prime}(1),\quad 1<\alpha<2,
\end{equation}It follows from (2.28)-(2.30) that the conditions of
Lemma 2.2 are fulfilled. Hence the proposition is true when
$1<\alpha<2.$ Let us consider the case when $0<\alpha<1$. As in the
previous case the function $\nu(x)$ is defined by relation (2.29).
We introduce the functions
\begin{equation}
k_{0}(y)=\int_{-\infty}^{y}\nu^{\prime}(t)dt\,y+\int_{y}^{0}\nu^{\prime}(t)tdt,
\quad y<0,\end{equation}
\begin{equation}
k_{0}(y)=-\int_{y}^{\infty}\nu^{\prime}(t)dt\,y-\int_{0}^{y}\nu^{\prime}(t)tdt,
\quad y>0,\end{equation} and  \begin{equation}
k(y)=k_{0}(y)+\gamma\mathrm{p}_{0}(y)\quad
0<\alpha<1.\end{equation}In view of (2.28) and (2.33),(2.34) the
conditions of Lemma
2.1 are fulfilled. Hence the proposition is proved.\\
\textbf{Corollary 2.1.} \emph{If condition $(2.28)$ is fulfilled then}
\begin{equation}
k_{0}(y){\geq}0,\quad
-\infty<y<\infty, \quad 1<\alpha<2,\end{equation}
\begin{equation}
k_{0}(y){\leq}0,\quad -\infty<y<\infty, \quad
0<\alpha<1.\end{equation} Let us consider the important case when
$\alpha=1.$\\
\textbf{Proposition 2.2.} \emph{The generator $L$ of the Levy
process $X_{t}$ admits the convolution type representation
$(2.1),(2.2)$ if there exist such $C>0$ and $m>0$ that}
\begin{equation}\nu^{\prime}(y){\leq}C|y|^{-2}e^{-m|y|}.\end{equation}
\emph{Proof.} Using formulas (2.29)-(2.32) we see that the
conditions of Lemma 2.2 are fulfilled. The proposition is proved.\\
\textbf{Example 2.1.} The stable processes.\\
For the stable processes we have $A=0,\quad
\gamma=\overline{\gamma}$ and
\begin{equation}
\nu^{\prime}(y)=|y|^{-\alpha-1}(C_{1}1_{y<0}+C_{2}1_{y>0}),
\end{equation}where $C_{1}>0\quad C_{2}>0.$ Hence the function $\nu(y)$ has the form
\begin{equation}\nu(y)=\frac{1}{\alpha}|y|^{-\alpha}(C_{1}1_{y<0}-C_{2}1_{y>0}).
\end{equation}Let us introduce the functions
\begin{equation}k_{0}(y)=\frac{1}{\alpha(\alpha-1)}|y|^{1-\alpha}(C_{1}1_{y<0}+C_{2}1_{y>0}),
\end{equation} where $0<\alpha<2,\quad \alpha{\ne}1.$ When $\alpha
=1$ we have \begin{equation}k_{0}(y)=
-\mathrm{log}|y|\,(C_{1}1_{y<0}+C_{2}1_{y>0}).\end{equation} It
means that the conditions of Theorem 2.1 are fulfilled. Hence the
generator $L$ for the stable processes admits the convolution type
representation (2.1),(2.2).\\
\textbf{Proposition 2.3.} \emph{The kernel $k(y)$ of the operator
$S$ in representation $(2.1)$ for the stable processes has  form
$(2.31)$, when $1{\leq}\alpha<2$, and has form $(2.35)$ when $0<\alpha<1$.} \\
\textbf{Example 2.2.} The variance damped Levy processes .\\
For the variance damped Levy processes we have $A=0,\quad
\gamma=\overline{\gamma}$ and \begin{equation}
\nu^{\prime}(y)=C_{1}e^{-\lambda_{1}
|y|}|y|^{-\alpha-1}1_{y<0}+C_{2}e^{-\lambda_{2}|y|}y^{-\alpha-1}1_{y>0},
\end{equation}where
$C_{1}>0\quad C_{2}>0,\quad \lambda_{1}>0,\quad \lambda_{2}>0 \quad
0<\alpha<2.$ It follows from (2.43) that the conditions of
Proposition 2.1 are fulfilled when $\alpha{\ne}1.$  If $\alpha=1$
the conditions of Proposition 2.2 are fulfilled. Hence the generator
$L$ for the variance damped Levy processes admits the convolution
type representation (2.1),(2.2) and the kernel $k(y)$ is defined by
formulas (2.30),(2.31), when $1{\leq}\alpha<2$,and by formula
(2.35) when  $0<\alpha<1$.\\
\textbf{Example 2.3.} The variance Gamma process .\\
For the variance Gamma process we have $A=0,\quad
\gamma=\overline{\gamma}$ and
\begin{equation}
\nu^{\prime}(y)=C_{1}e^{-G
|y|}|y|^{-1}1_{y<0}+C_{2}e^{-M|y|}y^{-1}1_{y>0},
\end{equation}where $C_{1}>0\quad C_{2}>0,\quad G>0,\quad
M>0.$
 It follows
from (2.44) that the conditions of Proposition 2.2 are fulfilled and
the generator $L$ of variance Gamma process admits the convolution
type representation (2.1),(2.2).
 The kernel $k(y)$ is defined by relations
(2.33) and (2.34).\\
\textbf{Example 2.4.} The normal inverse Gaussian process .\\
In the case of the normal inverse Gaussian process we have
$A=0,\quad \gamma=\overline{\gamma}$ and
\begin{equation}
\nu^{\prime}(y)=Ce^{\beta y}K_{1}(|y|)|y|^{-1},\quad C>0,\quad
-1{\leq}\beta{\leq}1,\end{equation}where $K_{\lambda}(x)$ denotes
the modified Bessel function of the third kind with index the
$\lambda$. Using equalities
\begin{equation}|K_{1}(|x|)|{\leq}Me^{-|x|}/|x|,\quad M>0,\quad
0<x_{0}{\leq}|x|,\end{equation}
\begin{equation}|K_{1}(|x|)x|{\leq}M,\quad
0{\leq}|x|{\leq}x_{0}\end{equation} we see that the conditions of
Proposition 2.2 are fulfilled. Hence the corresponding generator $L$
admits the convolution type representation (2.1), (2.2) and
 the kernel $k(y)$
is defined by  relations
(2.33) and (2.34).\\
\textbf{Example 2.5.} The Meixner process .\\
For the Meixner process we have \begin{equation}
\nu^{\prime}(y)=C\frac{\mathrm{exp}{\beta}x}{x\,\mathrm{sinh}{\pi}x}
,\end{equation}where $C>0,\quad -\pi<\beta<\pi.$ The conditions of
Proposition 2.2 are fulfilled. Hence the corresponding generator $L$
admits the convolution type representation (2.1),(2.2) and the
kernel $k(y)$ is defined by relations
(2.33),(2.34).\\
\textbf{Remark 2.1.} Examples 2.1-2.5 are used in the finance
problems [24].\\
\textbf{Example 2.6.} Compound Poisson process.\\
We consider the case when $A=0, \quad \gamma=0$ and
\begin{equation}
M=\int_{-\infty}^{\infty}\nu^{\prime}(y)dy<\infty.\end{equation}
 Using formulas (2.1) and (2.2) we deduce that the corresponding
 generator
$L$ has the following convolution form \begin{equation}
Lf=-Mf(x)+\int_{-\infty}^{\infty}\nu^{\prime}(y-x)f(y)dy.\end{equation}
\section{Potential}
The operator \begin{equation}
Qf=\int_{0}^{\infty}(P_{t}f)dt \end{equation}
is called \emph{potential} of the semigroup $P_{t}$. The generator $L$ and
the potential $Q$ are ( in general) unbounded operators. Therefore the operators
$L$ and $Q$ are defined not in the whole space $L^{2}(-\infty.\infty)$ but only
in the subsets $D_{L}$ and $D_{Q}$ respectively. We use the following property
of the potential $Q$ (see[23]).\\
\textbf{Proposition 3.1.} \emph{If $f=Q\,g,\quad (g{\in}D_{Q})$ then
$f{\in}D_{L}$ and }
\begin{equation}
-Lf=g.\end{equation}
\textbf{Example 3.1.} Compound Poisson process.\\
Let the generator $L$ has form (2.50) where \begin{equation}
M=\int_{-\infty}^{\infty}\nu^{\prime}(x)dx<\infty,\quad
\int_{-\infty}^{\infty}[\nu^{\prime}(x)]^{2}dx<\infty.
\end{equation} We introduce the functions \begin{equation}
K(u)=-\frac{1}{M\sqrt{2\pi}}\int_{-\infty}^{\infty}\nu^{\prime}(x)e^{-iux}dx,\end{equation}
\begin{equation}N(u)=\frac{K(u)}{1-\sqrt{2\pi}K(u)}.\end{equation}Let us note that
\begin{equation}
|K(u)|<\frac{1}{\sqrt{2\pi}},\quad u{\ne}0;\quad K(0)=-\frac{1}{\sqrt{2\pi}}.\end{equation}
It means that $N(u){\in}L^{2}(-\infty,\infty).$ Hence the function
\begin{equation}n(x)=-\frac{1}{\sqrt{2\pi}}\int_{-\infty}^{\infty}N(u)e^{-iux}du\end{equation}
belongs to $L^{2}(-\infty,\infty)$ as well. It follows from
(2.50),(3.2)and (3.7) that the corresponding potential $Q$ has the
form (see [23], Ch.11)\begin{equation}
Qf=\frac{1}{M}[f(x)+\int_{-\infty}^{\infty}f(y)n(x-y)dy].\end{equation}
\textbf{Proposition 3.2.} \emph{Let  conditions $(3.3)$  be
fulfilled. Then
the operators $L$ and $Q$ are bounded in the space}  $L^{2}(-\infty,\infty).$\\
Now we shall give an example when the kernel $n(x)$ can be written in an explicit form.\\
\textbf{Example 3.2.} We consider the case when \begin{equation}
{\nu}^{\prime}(x)=e^{-|x|},\quad -\infty<x<\infty.\end{equation}
In this case $M=2$ and the operator $L$ takes the form
\begin{equation}
Lf=-2f(x)+\int_{-\infty}^{\infty}f(y)e^{-|x-y|}dy.\end{equation}
Formulas (3.4)-(3.7) imply that
\begin{equation}
Qf=\frac{1}{2}f(x)-\frac{1}{4\sqrt{2}}\int_{-\infty}^{\infty}f(y)e^{-|x-y|\sqrt{2}}dy.\end{equation}
\section{Truncated generators and quasi-potentials}
Let us denote by $\Delta$ the set of segments $[a_{k},b_{k}]$ such that\\
 $a_{1}<b_{1}<a_{2}<b_{2}<...<a_{n}<b_{n},\quad 1{\leq}k{\leq}n.$ By $C_{\Delta}$
 we denote the set of functions $g(x)$ on $L^{2}(\Delta)$ such that
\begin{equation}g(a_{k})=g(b_{k})=g^{\prime}(a_{k})=g^{\prime}(b_{k})=0,\quad 1{\leq}k{\leq}n,
\quad g^{\prime\prime}(x){\in}L^{p}(\Delta),\quad p>1.
\end{equation}We introduce the operator $P_{\Delta}$ by relation
$P_{\Delta}f(x)=f(x)$ if $x{\in}\Delta$ and $P_{\Delta}f(x)=0$ if $x{\notin}\Delta.$\\
\textbf{Definition 4.1.} The operator
\begin{equation}
L_{\Delta}=P_{\Delta}LP_{\Delta} \end{equation}
is called a \emph{ truncated generator}.\\
\textbf{Definition 4.2.} The operator $B$ with the definition domain
dense in $L^{p}(\Delta)$  is called \emph{a quasi-potential} if the
functions $f=Bg$ belong to definition domain of $L_{\Delta}$ and
\begin{equation} -L_{\Delta}f=g.\end{equation} It follows from (4.3)
that
\begin{equation} -P_{\Delta}Lf=g,\quad ( f=Bg).\end{equation}
\textbf{Remark 4.1.} In a number of cases (see the next section) we need relation (4.4).
In these cases we can use the quasi-potential $B$, which is often simpler
than the corresponding potential $Q$.\\
\textbf{Remark 4.2.} The operators of type (4.2) are investigated in book ([22],Ch.2).\\
From relation (4.3) we deduce that
\begin{equation}
Bg{\ne}0,\quad if\quad g{\ne}0,\quad
g{\in}L^{p}(\Delta).\end{equation}
 \textbf{Definition 4.3.} We call
the operator $B$  a  \emph{regular} if the following
conditions are fulfilled.\\
1). The operator $B$ is compact and has the form \begin{equation}
Bf=\int_{\Delta}\Phi(x,y)f(y)dy,\quad f(y){\in}L^{p}(\Delta), \quad
p{\geq}1,\end{equation} where the function $\Phi(x,y)$ can have a
discontinuity only when $x=y.$ \\
2). There exists a function $\phi(x)$ such that
\begin{equation}|\Phi(x,y)|{\leq}\phi(x-y),\end{equation}
\begin{equation}\int_{-R}^{R}\phi(x)dx<\infty \quad if \quad
0<R<\infty.\end{equation} 3).\begin{equation}\Phi(x,y){\geq}0,\quad
x,y{\in}\Delta,\end{equation}
\begin{equation}\Phi(a_{k},y)=\Phi(b_{k},y)=0,\quad 1{\leq}k{\leq}n.\end{equation}
4). Relation (4.5) is valid.\\
\textbf{Remark 4.3.} In view of condition (4.7) the regular operator $B$
 is bounded  in the spaces
$L^{p}(\Delta),\quad 1{\leq}p{\leq}\infty$ (see[22],p.24).\\
\textbf{Remark 4.4.} If the quasi-potential $B$ is regular , then the corresponding
truncated generator $L_{\Delta}$ has a discrete spectrum. \\
Further we prove that for a broad class of Levy processes the
corresponding quasi-potentials $B$ are regular.\\
 \textbf{Example 4.1.} We consider the case when
\begin{equation} \phi(x)=M|x|^{-\varkappa},\quad
0<\varkappa<1.\end{equation} \textbf{Proposition 4.1.} \emph{Let
condition $(4.11)$ be true and let the corresponding regular
operator $B$ have an eigenfunction $f(x)$
 with an eigenvalue $\lambda{\ne}0$. Then the function $f(x)$ is continuous.}\\
 \emph{Proof.} According to  Definition 4.3 there exists an integer $N(\varkappa)$
 such that the kernel $\Phi_{N}(x,t)$ of the operator
\begin{equation}
B^{N}f=\int_{\Delta}\Phi_{N}(x,y)f(y)dy,\quad f(y){\in}L^{p}(\Delta)\end{equation}
is continuous. Hence the function $f(x)$ is continuous. The proposition is proved.

\section{The Probability of the Levy process remaining within the given domain}
In many theoretical and applied problems it is important to estimate the quantity
\begin{equation}p(t,\Delta)=P\{X_{\tau}{\in}\Delta\},\quad 0{\leq}\tau{\leq}t,\end{equation}
i.e. the the probability that a sample of the process $X_{\tau}$ remains inside $\Delta$
for $0{\leq}\tau{\leq}t$ (ruin problem).\\
To derive the integro-differential equations corresponding to Levy
processes we use the argumentation by Kac [11]  and our own
argumentation (see [20]-[22]). Now we get rid
 of the requirement for the process to be stable.\\
 Let us consider the Levy process $X_{t}$ with the continuous density
 (see (1.1):
\begin{equation}
\rho(x,t)=\frac{1}{2\pi}\int_{-\infty}^{\infty}e^{-ixz}\mu(z,t)dz,\quad t>0\end{equation}
Now we introduce the sequence of functions
\begin{equation}
Q_{n+1}(x,t)=\int_{0}^{t}\int_{-\infty}^{\infty}Q_{0}(x-\xi,t-\tau)V(\xi)Q_{n}(\xi,\tau)d{\xi}d\tau,
\end{equation} where the function $V(x)$ is defined by relations
 $V(x)=1$ when $x{\notin}\Delta$ and $V(x)=0$ when $x{\in}\Delta.$ We
 use the notation
 \begin{equation}Q_{0}(x,t)=\rho(x,t).\end{equation}
For Levy processes the following relation \begin{equation}
Q_{0}(x,t)=\int_{-\infty}^{\infty}Q_{0}(x-\xi,t-\tau)Q_{0}(\xi,\tau)d{\xi}\end{equation}
is true. Using (5.3) and (5.5) we have \begin{equation}
0{\leq}Q_{n}(x,t){\leq}t^{n}Q_{0}(x,t)/n!.\end{equation} Hence the
series \begin{equation}
Q(x,t,u)=\sum_{n-0}^{\infty}(-1)^{n}u^{n}Q_{n}(x,t) \end{equation}
converges. The probabilistic meaning of $Q(x,t,u)$ is defined by the
relation (see [12],Ch.4)
\begin{equation}E\{\mathrm{exp}[-u\int_{0}^{t}V(X_{\tau})d\tau],c_{1}<X_{t}<c_{2}\}=
\int_{c_{1}}^{c_{2}}Q(x,t,u)dx.\end{equation}The inequality $V(x){\geq}0$ and
relation (5.8) imply that the function $Q(x,t,u)$ monotonically decreases with respect
 to the variable "u" and the formulas \begin{equation}
 0{\leq}Q(x,t,u){\leq}Q(x,t,0)=Q_{0}(x,t)=\rho(x,t)\end{equation}
 are true. In view of (5.2) and (5.9) the Laplace transform
\begin{equation}
\psi(x,s,u)=\int_{0}^{\infty}e^{-st}Q(x,t,u)dt,\quad s>0.\end{equation}
has the meaning.
According to (5.3) the function $Q(x,t,u)$ is the solution of the equation
\begin{equation}
Q(x,t,u)+u\int_{0}^{t}\int_{-\infty}^{\infty}\rho(x-\xi,t-\tau)V(\xi)Q(\xi,\tau,u)d{\xi}d\tau=
\rho(x,t)\end{equation}
Taking from both parts of (5.11) the Laplace transform and bearing in mind  (5.10) we obtain
\begin{equation}
\psi(x,s,u)+u\int_{-\infty}^{\infty}V(\xi)R(x-\xi,s)\psi(\xi,s,u)d\xi=R(x,s),
\end{equation} where
\begin{equation}
R(x,s)=\int_{0}^{\infty}e^{-st}\rho(x,t)dt.\end{equation}Multiplying both parts
of relation
(5.12) by $\mathrm{exp}(ixp)$ and integrating them with respect to
$x \quad (-\infty<x<\infty)$ we have \begin{equation}
\int_{-\infty}^{\infty}\psi(x,s,u)e^{ixp}[s+\lambda(p)+uV(x)]dx=1.\end{equation}
Here we use relations (1.1), (5.2) and (5.13). Now we introduce the function
\begin{equation}h(p)=\frac{1}{2\pi}\int_{\Delta}e^{-ixp}f(x)dx,\end{equation}
where the function $f(x)$ belongs to $C_{\Delta}$. Multiplying both parts  of (5.14)
by $h(p)$ and integrating them with respect to $p\quad  (-\infty<p<\infty)$ we deduce
the equality
\begin{equation}
\int_{-\infty}^{\infty}\int_{-\infty}^{\infty}\psi(x,s,u)e^{ixp}[s+\lambda(p)]h(p)dxdp=f(0).
\end{equation} We have used the relations \begin{equation}
V(x)f(x)=0,\quad -\infty<x<\infty,\end{equation}
\begin{equation}
\frac{1}{2\pi}\mathrm{lim}\int_{-N}^{N}\int_{\Delta}e^{-ixp}f(x)dxdp=f(0),\quad N{\to}\infty.
\end{equation}
Since the function $Q(x,t,u)$ monotonically decreases with respect
to $"u"$
 this is also true for the function $\psi(x,s,u)$ according to (5.10). Hence there exists the limit
 \begin{equation}
 \psi(x,s)=\mathrm{\lim}\psi(x,s,u),\quad u{\to}\infty,\end{equation}where
 \begin{equation}
 \psi(x,s)=0,\quad x{\notin}\Delta.\end{equation}
 The probabilistic meaning of $\psi(x,s)$
 follows from the equality \begin{equation}
 \int_{0}^{\infty}e^{-st}p(t,\Delta)dt=\int_{\Delta}\psi(x,s)dx.\end{equation}
 Using the properties of the Fourier transform and  conditions (5.19) , (5.20)
  we deduce from (5.16) the following assertion.\\
  \textbf{Theorem 5.1.} \emph{Let the considered Levy process
   have the continuous density.
 Then the relation  \begin{equation}
  ((sI-L_{\Delta})f,\psi(x,s))_{\Delta}=f(0)\end{equation} is true.}\\
  \textbf{Remark 5.1} For symmetric stable processes relation (5.22)
  was deduced by M.Kac [11].\\
  \textbf{Remark 5.2} As it is known ,the stable processes, the
  variance damped Levy processes , the variance gamma processes, the
  normal inverse Gaussian process, the Meixner process have
  continuous densities (see([24],[31]).\\
 \textbf{Remark 5.3.} So we have obtained the formula (5.21) for
   Laplace transform of $p(t, \Delta)$
in terms of $\psi(x,s)$. The double Laplace transform of
$p(t,\Delta)$ was obtained
 by G.Baxter and M.D.Donsker [3] for the case when $\Delta=(-\infty,a].$\\
    We express
  the important function $\psi(x,s)$   with the help of the
  quasi-potential $B$.\\
\textbf{Theorem 5.2.} \emph{Let the considered Levy process have
  the continuous density and let the
  quasi-potential $B$ be regular .
  Then in the space $L^{p}(\Delta)\quad (p>1)$ there is one and only one function
  \begin{equation}
  \psi(x,s)=(I+sB^{\star})^{-1}\Phi(0,x) , \quad 0{\leq}s<s_{0},\end{equation}
  which satisfies relation $(5.22)$.}\\
  \emph{Proof.} In view of (4.4) we have   \begin{equation}
  -BL_{\Delta}f=f,\quad f{\in}C_{\Delta}. \end{equation}Relations (5.23) and (5.24)
  imply that \begin{equation}
((sI-L_{\Delta})f,\psi(x,s))_{\Delta}=-((I+sB)L_{\Delta}f,\psi)_{\Delta}=-(L_{\Delta}f,\Phi(0,x))_{\Delta}.\end{equation}
Since $\Phi(0,x)=B^{\star}\delta(x),$  ($\delta(x)$  is the Dirac
function) then according to (5.23) and (5.25) relation (5.22) is
true.\\
 Let us suppose that in $L(\Delta)$ there is another function
$\psi_{1}(x,s)$ satisfying (5.22). Then the equality
\begin{equation} ((sI-L_{\Delta})f,\phi(x,s))_{\Delta}=0,\quad
\phi=\psi-\psi_{1} \end{equation} is valid. We write  relation
(5.26) in the form
\begin{equation}
(L_{\Delta}f,(I+sB^{\star})\phi)_{\Delta}=0.\end{equation} Due to
(4.4) the range of $L_{\Delta}$ is dense in $L^{p}(\Delta).$ Hence
in view of (5.27) we have
$\phi=0.$ The theorem is proved.\\
The analytical apparatus for the construction and investigation of
the function $\psi(x,s)$ is based on  relation (5.22) and properties
of the quasi-potential $B$. In the following three sections
 we shall investigate the properties of the operator $B$.
\section{Non-negativity of the kernel $\Phi(x,y)$}
 In this section we deduce the following important property of the kernel
$\Phi(x,y)$.\\
 \textbf{Proposition 6.1.} \emph{Let the density
$\rho(x,t)$ of Levy process $X_{t}$ be continuous   $(t>0)$ and let
the corresponding quasi-potential $B$ satisfy conditions
$(4.6)-(4.8)$  of Definition $4.3.$
 Then the  kernel $\Phi(x,y)$
is non-negative i.e.} \begin{equation}\Phi(x,y){\geq}0.\end{equation}
\emph{Proof.} In view of (5.9) and (5.10) we have $ \psi(x,s,u){\geq}0$. Relation
(5.19)implies that  $ \psi(x,s){\geq}0$. Now it follows from (5.23) that
\begin{equation}\Phi(0,x)=\psi(x,0){\geq}0.\end{equation}
Let us consider the domains $\Delta_{1}$ and $\Delta_{2}$ which are
connected by relation $\Delta_{2}=\Delta_{1}+\delta$. We denote the
corresponding truncated generators  by $L_{\Delta_{1}}$  and
$L_{\Delta_{2}}$, we denote the corresponding quasi-potentials by
$B_{1}$ and $B_{2}$ and the corresponding kernels by $\Phi_{1}(x,y)$
and $\Phi_{2}(x,y)$. We introduce the unitary operator
\begin {equation} Uf=f(x-\delta),
\end{equation} which maps the space $L^{2}(\Delta_{2})$ onto
$L^{2}(\Delta_{1})$. At the beginning we suppose that the conditions
of Theorem 2.1 are fulfilled. Using formulas (2.1) and (2.2) we
deduce that
\begin{equation} L_{\Delta_{2}}=U^{-1}L_{\Delta_{1}}U.\end{equation}
Hence the equality
\begin{equation}
B_{2}=U^{-1}B_{1}U\end{equation} is valid. The last equality can be
written in the terms of the kernels \begin{equation}
\Phi_{2}(x,y)=\Phi_{1}(x+\delta,y+\delta).\end{equation}According to
(6.2) and (6.6) we have
\begin{equation}\Phi_{1}(\delta,y+\delta){\geq}0.\end{equation} As
$\delta$  is an arbitrary real number, relation (6.1) follows
directly from (6.6).
 We remark that an arbitrary generator  $L$ can be approximated
 by the operators of form
(2.1) (see[23],Ch.2). Hence the proposition is proved.\\
In view of (4.1), (4.5) and relation $Bf{\in}C_{\Delta}$  the following assertion is true.\\
\textbf{Proposition 6.2.} \emph{Let the quasi-potential $B$ satisfy
the conditions of Proposition $6.2$, Then the equalities}
\begin{equation}
\Phi(a_{k},y)=\Phi(b_{k},y)=0 \quad 1{\leq}k{\leq}n \end{equation}
\emph{are valid}.\\
\section{Sectorial operators}
1. We introduce the following notions.\\
\textbf{Definition 7.1.} The bounded operator $B$ in the space $L^{2}(\Delta)$ is called
\emph{sectorial} if \begin{equation}
(Bf,f){\ne}0,\quad f{\ne}0 \end{equation}
and
\begin{equation}
-\frac{\pi}{2}\beta{\leq}\mathrm{arg}(Bf,f){\leq}\frac{\pi}{2}\beta, \quad 0<\beta{\leq}1.
\end{equation}
It is easy to see that the following assertions are true.\\
\textbf{Proposition 7.1.} \emph{ Let the operator $B$ be sectorial.
Then the operator $(I+sB)^{-1}$
is bounded when} $s{\geq}0$.\\
\textbf{Proposition 7.2.} \emph{Let the conditions of Theorem $5.2$
be fulfilled. If the operator $B$ is sectorial, then formula
$(5.23)$ is
valid for all }$s{\geq}0.$\\
In the present section we deduce the conditions under which the
quasi-potential $B$ is  sectorial. Let us consider the case when
\begin{equation} \int_{x}^{\infty}y{\nu}^{\prime}(y)dy<\infty,\quad (x>0),\end{equation}
\begin{equation} \int_{-\infty}^{x}|y|{\nu}^{\prime}(y)dy<\infty,\quad (x<0).\end{equation}
The corresponding kernel $k(x)$ of the operator $S$ (see(2.2) has the form
\begin{equation}
k(x)=\int_{x}^{\infty}(y-x){\nu}^{\prime}(y)dy<\infty,\quad (x>0),\end{equation}
\begin{equation}k(x)= \int_{-\infty}^{x}(x-y){\nu}^{\prime}(y)dy<\infty,\quad (x<0).\end{equation}
Using the inequality ${\nu}^{\prime}(y){\geq}0$ we obtain the following statement.\\
\textbf{Proposition 7.3.} \emph{Let conditions $(7.3)$ and $(7.4)$
be fulfilled . Then the kernel $k(x)$ is monotone on the half-axis
$(-\infty,0)$ and on the half-axis}
$(0,\infty)$.\\We shall use the following Pringsheim's result .\\
\textbf{Theorem 7.1.} (see[25], Ch.1) \emph{ Let $f(t)$ be
non-increasing function over $(0,\infty)$ and integrable on any
finite interval }$(0,\ell)$. \emph{If $f(t){\to}0$ when
$t{\to}\infty$, then for any positive $x$ we have}
\begin{equation}\frac{1}{2}[f(x+0)+f(x-0)]=\frac{2}{\pi}\int_{+0}^{\infty}\mathrm{cos}xu
[\int_{0}^{\infty}f(t)\mathrm{cos}tudt]du,\end{equation}
\begin{equation}\frac{1}{2}[f(x+0)+f(x-0)]=\frac{2}{\pi}\int_{0}^{\infty}\mathrm{sin}xu
[\int_{0}^{\infty}f(t)\mathrm{sin}tudt]du.\end{equation}
It follows from (7.3)-(7.6) that
\begin{equation} k(x){\to}0 \quad and \quad k^{\prime}(x){\to}0,\quad
when\quad x{\to}{\pm}\infty.\end{equation}
We suppose in addition that
\begin{equation} xk(x){\to}0 \quad and \quad x^{2}k^{\prime}(x){\to}0,\quad
when\quad x{\to}{\pm}0.\end{equation}
Using the integration by parts we deduce the assertion.\\
\textbf{Proposition 7.4.} \emph{Let conditions $(7.3)$ , $(7.4)$ and
$(7.9)$, $(7.10)$ be fulfilled. Then the relation}
\begin{equation}
\int_{-\infty}^{\infty}k(t)\mathrm{cos}xtdt=\int_{-\infty}^{\infty}{\nu}^{\prime}(t)
\frac{1-\mathrm{cos}xt}{x^{2}}dt \end{equation}
\emph{is true.}\\
Relation (7.11) implies that
\begin{equation}\int_{-\infty}^{\infty}k(t)\mathrm{cos}xtdt>0.\end{equation}
The kernel $k(x)$ of the operator $S$ admits the representation
\begin{equation} k(x)=\int_{-\infty}^{\infty}m(t)e^{ixt}dt.\end{equation}
In view of (7.12) we have
\begin{equation} \mathrm{Re}[m(u)]>0.\end{equation}
Due to (7.13) and (7.14) the relation
\begin{equation}
(Sf,f)=\int_{-\infty}^{\infty}m(u)|\int_{\Delta}f(t)e^{iut}dt|^{2}du\end{equation}
is valid. Hence we have
\begin{equation}
-\frac{\pi}{2}{\leq}\mathrm{arg}(Sf,f){\leq}\frac{\pi}{2},
\quad f(t){\in}L^{2}(\Delta).
\end{equation}
\textbf{Proposition 7.5.} \emph{Let conditions $(7.3)$ , $(7.4)$ and
$(7.9)$,  $(7.10)$ be fulfilled.
Then the corresponding operator $B$ is sectorial.}\\
\emph{Proof.} Let the function $g(x)$ satisfies conditions (4.1). Then the relation
\begin{equation}
(-Lg,g)=(Sg^{\prime},g^{\prime}) \end{equation}
holds. Equalities (4.3) and (7.17) imply that
\begin{equation}(f,Bf)=(Sg^{\prime},g^{\prime}) ,\quad g=Bf.\end{equation}
Inequality (7.1) follows
from relations (7.14) and (7.18).
 Relations (7.16) and (7.18) imply the proposition.\\
\textbf{Remark 7.1.} The variance damped processes ( Example 2.2.)
the normal inverse Gaussian process (Example 2.4.), the Meixner
process (Example  2.5.) satisfy the conditions of Proposition 7.5.
Hence the corresponding operators $B$
 are  sectorial.\\
2. Now  we introduce the notion of the strongly sectorial operators.\\
 \textbf{Definition 7.2.} The sectorial operator $B$ is called
 \emph{a strongly sectorial}
 if for some $\beta<1$  relation
(7.2)
is valid.\\
\textbf{Proposition 7.6.} \emph{Let the following conditions be fulfilled.\\
1). Relations $(7.3)$ , $(7.4)$ and $(7.9)$, $(7.10)$ are valid.}\\
2).\emph{ For some $m>0$  the inequality}
\begin{equation}
\frac{m}{|x|^{2}}{\leq}\nu^{\prime}(x),\quad |x|{\leq}1 \end{equation}
\emph{is true.}\\
3.\begin{equation}\int_{-\infty}^{\infty}k(t)dt<\infty.\end{equation}
\emph{Then the corresponding operator $B$ is strongly sectorial .}\\
\emph{Proof.}  As it is known (see [25],Ch.1) the inequality
\begin{equation}|\int_{-\infty}^{\infty}k(t)\mathrm{sin}xtdt|{\leq}\frac{M}{|x|},
\quad M>0,\quad |t|{\geq}1 \end{equation} is valid. From formulas (7.11) and (7.26)
 we conclude that
\begin{equation}
\int_{-\infty}^{\infty}k(t)\mathrm{cos}xtdt{\geq}\int_{-1/x}^{1/x}{\nu}^{\prime}(t)
\frac{1-\mathrm{cos}xt}{x^{2}}dt{\geq}\frac{N}{|x|},\quad N>0,\quad |x|{\geq}1.\end{equation}
It follows from (7.28) and (7.29) that
\begin{equation}
-\frac{\pi}{2}\beta{\leq}\mathrm{arg}(Sf,f){\leq}\frac{\pi}{2}\beta, \quad 0<\beta<1.
\end{equation}
Hence according to (7.18) relation (7.25) is valid .
The proposition is proved.\\
\textbf{Remark 7.2.} The variance damped processes (Example 2.2,
$\alpha{\geq}1$), the normal inverse Gaussian process (Example
2.4.), the Meixner process (Example  2.5.) satisfy the conditions of
Proposition 7.6. Hence the corresponding operators $B$
 are  strongly sectorial.\\
 \textbf{Proposition 7.7.} \emph{Let conditions $(7.3), (7.4)$ and $(7.9)$, $(7.10)$ be fulfilled.
If the operator $S$ has the form}
\begin{equation}
Sf=Af+\int_{\Delta}k(x-t)f(t)dt,\quad A>0.\end{equation}
\emph{Then the corresponding operator $B$ is strongly sectorial.}\\
\emph{Proof.} It is easy to see that for some $\beta<1$  relation
(7.23)  is true. According to relation (7.18) the corresponding
operator $B$ is strongly sectorial.
\section{Quasi-potential $B$, structure and properties}
Let us begin  with the symmetric segment $\Delta=[-c,c]$.\\
\textbf{Theorem 8.1.} (see[22],p.140) \emph{Let the following
 conditions be fulfilled\\
1. There exist the functions $N_{k}(x){\in}L^{p}(-c,c),\quad p>1$
 which satisfy the equations
\begin{equation}
SN_{k}=x^{k-1},\quad k=1,2.\end{equation} 2.\begin{equation}
r=\int_{-c}^{c}N_{1}(x)dx{\ne}0 \end{equation} } \emph{Then the
corresponding operator $B$ has the form
\begin{equation}Bf=\int_{-c}^{c}\Phi(x,y,c)f(y)dy\end{equation}
where}
\begin{equation}
\Phi(x,y,c)=\frac{1}{2}\int_{x+y}^{2c-|x-y|}q[(s+x-y)/2,(s-x+y)/2]ds,\end{equation}
\begin{equation}q(x,y)=[N_{1}(-y)N_{2}(x)-N_{2}(-y)N_{1}(x)]/r.\end{equation}
It follows from (8.4) and (8.5) that \begin{equation}
\Phi({\pm}c,y)=\Phi(x,{\pm}c)=0.\end{equation} Here we use the
following relation\\
$q[(s+x-y)/2,(s-x+y)/2]=$
\begin{equation}[N_{1}((x-y-s)/2)N_{2}((s+x-y)/2)-N_{2}((x-y-s)/2)N_{2}((s+x-y)/2)]/r.
\end{equation}Thus
\begin{equation}q[(s+x-y)/2,(s-x+y)/2]=-q[(-s+x-y)/2,(-s-x+y)/2].\end{equation}
From formulas (8.4) and (8.5) we deduce the following statement.\\
\textbf{Proposition 8.1.} \emph{Let the conditions of Theorem $8.1$
be fulfilled.  There exists a function $\phi(x)$ such that}
\begin{equation}|\Phi(x,y,c)|{\leq}\phi(x-y),\end{equation}
\begin{equation}\int_{-R}^{R}\phi(x)dx<\infty \quad if \quad
0<R<\infty.\end{equation} \emph{Proof.} Relation (8.4) can be
written in the form
\begin{equation}\Phi(x,y,c)=\int_{x}^{c+(x-y-|x-y|)/2}q(t,t-x+y)dt.\end{equation}
By relations
\begin{equation}N_{k}(x)=0,\quad x{\notin}[-c,c], \quad k=1,2
\end{equation} we extend the functions $N_{k}(x)$ from the segment
$[-c,c]$ to the segment $[-2c,2c].$ It follows from (8.11) and
(8.12) that inequality (8.9) is valid , if
\begin{equation}
\phi(x)=\int_{-c}^{c}[|N_{1}(t)N_{2}(t-x)|+|N_{2}(t)N_{1}(t-x)|]dt/|r|.\end{equation}
Equality (8.13) imply that $\phi(x){\in}L^{p}[-2c,2c].$ The
proposition is proved.\\
 It follows from Proposition 8.1 that the
operator $B$ is bounded in all the spaces $L^{p}(-c,c),\quad
p{\geq}1$.  We shall prove that the operator $B$ is compact.\\
\textbf{Proposition 8.2.} \emph{Let the conditions of Theorem $8.1$
be fulfilled. Then the operator $B$ is compact in all the spaces}
$L^{p}(-c,c),\quad p{\geq}1.$\\
\emph{Proof.} Let us consider the operator $B^{\star}$ in the space
$L^{q}(-c,c),\quad 1/p+1/q=1.$ Using relation (8.3) we have
\begin{equation} B^{\star}f_{n}=\int_{-c}^{c}\Phi(y,x,c)f_{n}(y)dy
\end{equation} where the functions $f_{n}(x){\to}0$ in the weak sense.
Relation (8.14) can be represented in the following form
\begin{equation}B^{\star}f_{n}=\int_{-c}^{c}f_{n}(y)\int_{y}^{c+(y-x-|x-y|)/2}q(t,t-y+x)dtdy.
\end{equation}By interchanging the order of the integration in (8.15) we see
that $||B^{\star}f_{n}||{\to}0,$ i.e. the operator $B^{\star}$ is
compact. Hence the operator $B$ is compact too. The proposition is
proved.\\
Using formulas (8.5) and (8.11) we obtain the assertion.\\
\textbf{Proposition 8.3} \emph{Let the conditions of Theorem $8.1$
be fulfilled. If the functions $N_{1}(x)$ and $N_{2}(x)$ can have a
discontinuity only when $x={\pm}c$ then the function $\Phi(x,y,c)$
can have a discontinuity only when} $x=y.$\\
\textbf{Corollary 8.1.} \emph{Let the conditions of Proposition
$8.3$ be fulfilled. Then the eigenvectors of the corresponding
operator $B$ are continuous.}\\We use the following assertion
(see[22],p.73).\\
\textbf{Proposition 8.4.} \emph{Let the following conditions be fulfilled.\\
1).The kernel $k(x)$ of the operator $S$ has the form}
\begin{equation}
k(x)=\mathrm{log}\frac{A}{2|x|}+h(x),\quad -2c{\leq}x{\leq}2c,
\end{equation}
\emph{where} $ A>0,\quad A{\ne}c.$\\
2).
\begin{equation}\int_{-2c}^{2c}|h^{\prime}(u)|^{q}(2c-|u|)du<\infty,\quad
q>2.\end{equation} 3.\emph{ The equation}
\begin{equation}
Sf=\int_{-c}^{c}[\mathrm{log}\frac{A}{2|x-t|}+h(x-t)]f(t)dt=0
\end{equation}
\emph{has only the trivial solution in $L^{p}(-c,c)$, where $1/p+1/q=1.$\\
Then the equation}
\begin{equation}Sf=g,\quad g^{\prime}(x){\in}L^{p}(-c,c) \end{equation}
\emph{has one and only one solution in} $L^{p}(-c,c).$\\
\textbf{Corollary 8.2.} \emph{Let the conditions of Proposition
$8.4$ be valid. Then there exist  the functions $N_{1}(x)$ and
$N_{2}(x)$ which satisfy the equations}
\begin{equation}SN_{k}=x^{k-1},\quad N_{k}(x){\in}L^{p}(-c,c),\quad
k=1,2.\end{equation} \emph{ The functions $N_{1}(x)$ and $N_{2}(x)$
can
have a discontinuity only when} $x={\pm}c.$ \\
\textbf{Remark 8.1.} If conditions (7.13) and (7.14) are fulfilled
then according to (7.15) we have $(f,Sf){\ne}0$,\quad when
$||f||{\ne}0$. In particular the relation $(N_{1},SN_{1})=r{\ne}0$
is true.\\
\textbf{Remark 8.2.} In view of (6.4) and (6.5) Proposition 8.1 is
valid not only in the
case of the symmetric segment $[-c,c]$ but in the general case $[-a,b]$ too.\\
\section{Long time  behavior}
1. In order to investigate the asymptotic behavior of $p(t,\Delta)$ when
$t{\to}\infty$, we use the non-negativity of the kernel $\Phi(x,y)$. We apply the
following Krein-Rutman theorem (see [13], section 6).\\
\textbf{Theorem 9.1.} \emph{If a linear compact operator $B$ leaving invariant a cone
$K$,  has a point of the spectrum different from zero , then it has a positive eigenvalue $\lambda_{1}$
not less in modulus than any other eigenvalues $\lambda_{k},\quad (k>1)$.
 To this eigenvalue $\lambda_{1}$ corresponds
at least one eigenvector $g_{1}{\in}K, (Bg_{1}=\lambda_{1}g_{1})$ of
the operator $B$ and at least one eigenvector $h_{1}{\in}K^{\star},
(B^{\star}h_{1}=\lambda_{1}h_{1})$ of the operator $B^{\star}$.}\\
We remark that in our case the cone $K$ consists of non-negative
functions $f(x){\in}L^{p}(\Delta)$. Hence we have
\begin{equation}
g_{1}(x){\geq}0,\quad h_{1}(x){\geq}0.\end{equation}
 We introduce the
 following normalizing condition
\begin{equation}
 (g_{1},h_{1})=\int_{\Delta}g_{1}(x)h_{1}(x)dx=1,\end{equation}
Let the interval $\Delta_{1}$ and the point $x_{0}$ be such that
 \begin{equation}
 x_{0}{\in}\Delta_{1}{\in}\Delta.\end{equation} Together with quantity $p(t,\Delta)$
  we consider the expression
  \begin{equation}
  p(x_{0},\Delta_{1},t,\Delta)=P(\underset{0{\leq}\tau{\leq}t}{(X_{\tau}{\in}\Delta)}{\bigcap}(X_{t}{\in}\Delta_{1})),
 \end{equation}where $x_{0}=X_{0}$. If the relations $x_{0}=0,\quad \Delta_{1}=\Delta $ are true, then
 $p(x_{0},\Delta_{1},t,\Delta)=p(t,\Delta).$
   In this section we investigate the asymptotic
 behavior of  $p(x_{0},\Delta_{1},t,\Delta)$ and   $p(t,\Delta)$ when $t{\to}\infty$.\\
\textbf{Theorem 9.2.} \emph{Let the considered Levy process have the continuous density and let
 the corresponding quasi-potential $B$ be regular and  strongly
 sectorial. And  let the operator $B$ have a point of the
 spectrum different from zero.
 Then the asymptotic equality }\begin{equation}
 p(t,\Delta)=e^{-t/\lambda_{1}}[q(t)+o(1)],\quad t{\to}+\infty
 \end{equation}\emph{is true. The function $q(t)$ has the form}
 \begin{equation}
 q(t)=c_{1}+\sum_{k=2}^{m}c_{k}e^{it\nu_{k}}{\geq}0,\end{equation}
\emph{where $\nu_{k}$ are real}\\
 \emph{Proof.} The spectrum
$(\lambda_{k},\quad k>1)$ of the operator $B$ is situated
 in the sector
\begin{equation}
-\frac{\pi}{2}\beta{\leq}\mathrm{arg}z{\leq}\frac{\pi}{2}\beta,\quad 0{\leq}\beta<1,
\quad |z|{\leq}\lambda_{1}.
\end{equation}We introduce the domain $D_{\epsilon}$:
\begin{equation}
-\frac{\pi}{2}(\beta
+\epsilon){\leq}\mathrm{arg}z{\leq}\frac{\pi}{2}(\beta
+\epsilon),\quad |z-(1/2)\lambda_{1}|<(1/2)(\lambda_{1}
-r),\end{equation} where $0<\epsilon<1-\beta),\quad r<\lambda_{1}$.
If $z$ belongs to the domain $D_{\epsilon}$ then the relation
\begin{equation} \mathrm{Re}(1/z)>1/\lambda_{1}\end{equation}
holds. As the operator $B$ is compact only a finite number of
eigenvalues $\lambda_{k}, \quad 1<k{\leq}m$ of this operator does
not belong to the domain $D_{\epsilon}$. We denote  the boundary of
domain $D_{\epsilon}$
 by $\Gamma_{\epsilon}$. Without loss of generality we may assume that the points of spectrum
 $\lambda_{k}{\ne}0$ do not belong to $\Gamma_{\epsilon}$.
 Taking into account the equality
 \begin{equation}
 (\Phi(0,x),g_{1}(x))=\lambda_{1}g_{1}(0), \end{equation}
 we deduce from formulas (5.21) and (5.23) the relation
 \begin{equation}
 p(t,\Delta)=
 \sum_{k=1}^{m}\sum_{j=0}^{n_{k}}e^{-t/\lambda_{k}}t^{j}c_{k,j}+J,\end{equation}
where $n_{k}$ is the index of the eigenvalue
$\lambda_{k}$,
\begin{equation}
J=-\frac{1}{2i\pi}\int_{\Gamma}\frac{1}{z}e^{-t/z}((B^{\star} -zI)^{-1}\Phi(0,x),1)dz.
\end{equation}We note that \begin{equation}n_{1}=1.\end{equation}
Indeed, if $n_{1}>1$ then there exists such a function $f_{1}$ that
\begin{equation}
Bf_{1}=\lambda_{1}f_{1}+g_{1}.\end{equation}
In this case the relations
\begin{equation}
(Bf_{1},h_{1})=\lambda_{1}(f_{1},h_{1})+(g_{1},h_{1})=\lambda_{1}(f_{1},h_{1})
\end{equation}are true. Hence $(g_{1},h_{1})=0.$ The last relation contradicts
 condition (9.2). It proves equality (9.13).\\
 Relation (8.9) implies  that \begin{equation}
 \Phi(0,x){\in}L^{p}(\Delta).\end{equation}
We denote by $W(B)$ the numerical range of $B$. The closure of the
convex hull of $W(B)$ is situated in the sector (9.7). Hence the
estimation
\begin{equation} ||(B^{\star}-zI)^{-1}||_{p}{\leq}M/|z|,\quad z{\in}\Gamma_{\epsilon} \end{equation}
is true (see [26] for the Hilbert case $p=2$ and  [16],[28] for the
Banach space $p{\geq}1$).
 By $||B||_{p}$ we denote
the norm of the operator $B$ in the space $L^{p}(\Delta)$.\\ It
follows from estimation(9.17)
that the integral J exists.\\
Among the numbers $\lambda_{k}$ we choose  for which
$\mathrm{Re}(1/{\lambda}_{k}),\quad (1{\leq}k{\leq}m)$ has the
smallest value $\delta$. Among the obtained numbers we choose
$\mu_{k},\quad (1{\leq}k{\leq}\ell)$ the indexes $n_{k}$ of which
have the largest value $n$. We deduce from (9.10)-(9.12) that
\begin{equation}
 p(t,\Delta)=e^{-t{\delta}}t^{n}
 [\sum_{k=1}^{\ell}e^{-t/{\mu}_{k}}c_{k}+o(1)],\quad t{\to}\infty.\end{equation}
 We note that the function
 \begin{equation} Q(t)= \sum_{k=1}^{\ell}e^{it\mathrm{Im}(\mu_{k}^{-1})}c_{k} \end{equation}
 is almost periodic (see [14]). Hence in view of (9.18) and the inequality\\
  $p(t,\Delta)>0,\quad t{\geq}0$ the following relation
 \begin{equation} Q(t){\geq}0,\quad -\infty<t<\infty \end{equation}
 is valid.\\
 First we assume that at least one of the inequalities
 \begin{equation}\delta<{\lambda}_{1}^{-1},\quad n>1\end{equation}
 is true. Using  (9.21) and the inequality
 \begin{equation} \lambda_{1}{\geq}\lambda_{k},\quad k=2,3,...\end{equation}we have
\begin{equation}\mathrm{Im}{\mu}_{j}^{-1}{\ne}0,\quad 1{\leq}j{\leq}\ell.\end{equation}
It follows from (9.19) that
\begin{equation}
c_{j}=\mathrm{lim}\frac{1}{2T}\int_{-T}^{T}Q(t)e^{-it(\mathrm{Im}{\mu}_{j}^{-1})}dt,\quad
T{\to}\infty.\end{equation}In view of (9.20) we obtain the relations
\begin{equation}
|c_{j}|{\leq}\mathrm{lim}\frac{1}{2T}\int_{-T}^{T}Q(t)dt=0,\quad
T{\to}\infty,\end{equation}i.e. $c_{j}=0,\quad 1{\leq}j{\leq}\ell.$
This means that  relations (9.21) are not true. Hence the equalities
 \begin{equation}\delta={\lambda}_{1}^{-1},\quad n=1\end{equation}
 are true. From (9.18) and (9.19) we get the asymptotic equality
 \begin{equation}
 p(t,\Delta)=e^{-t/\lambda_{1}}[q(t)+o(1)]\quad t{\to}\infty,\end{equation}
 where the function $q(t)$ is defined by relation (9.6) and
 \begin{equation} c_{k}=\overline{g_{k}(0)}\int_{\Delta}h_{k}(x)dx,
 \quad \nu_{k}=\mathrm{Im}(\mu^{-1}).\end{equation}
 Here $g_{k}(x)$ are the eigenfunctions of the operator $B$ corresponding to
  the eigenvalues $\lambda_{k}$, and $h_{k}(x)$ are the eigenfunctions of the operator $B^{\star}$
   corresponding to
  the eigenvalues $\overline{\lambda_{k}}.$ The following conditions are fulfilled
  \begin{equation}
 (g_{k},h_{k})=\int_{\Delta}\overline{g_{k}(x)}h_{k}(x)dx=1,\end{equation}
\begin{equation}
 (g_{k},h_{\ell})=\int_{\Delta}\overline{g_{k}(x)}h_{\ell}(x)dx=0,\quad k{\ne}\ell.\end{equation}
  Using the almost periodicity of the function
  $q(t)$ we deduce from (9.27) the inequality
  \begin{equation}q(t){\geq}0.\end{equation} The theorem is proved.\\
 \textbf{Corollary 9.1.} \emph{ Let the conditions of Theorem 9.2 be fulfilled.
 Then all the eigenvalues $\lambda_{j}$ of $B$
belong  to the disk}
 \begin{equation}
|z-(1/2)\lambda_{1}|{\leq}(1/2).\end{equation}
\emph{All the eigenvalues $\lambda_{j}$ of $B$ which belong to the boundary of
 disc (9.32) have the indexes} $n_{j}=1.$\\
\textbf{Remark 9.1.} The exponential decay of the transition probability $P_{t}(x,B)$
 was proved by
P. Tuominen and R.L.Tweedie [29]. Theorem 9.2. gives the exponential
decay of $p(t,\Delta)$.
 These two results are independent.\\
Using formula (9.11) we obtain the following assertion.\\
 \textbf{Corollary 9.2.} \emph{Let the considered Levy process have the continuous density and let
 the corresponding quasi-potential $B$ be regular and  strongly
 sectorial. And  let the operator $B$ have no points of the
 spectrum different from zero.
 Then the  equality }\begin{equation}
 \mathrm{lim}[p(t,\Delta)e^{t/\lambda}]=0,\quad t{\to}+\infty \end{equation}
 \emph{ is true
 for any }$\lambda>0.$\\
2. Now we find the conditions under which the operator $B$ has a
point of the
 spectrum different from zero.\\
 We represent the corresponding operator $B$ in the form
 $B=B_{1}+iB_{2}$ where the operators $B_{1}$ and $B_{2}$ are self-adjoint.
 We assume that $B_{1}{\in}\Sigma_{p}$, i.e.
 \begin{equation}\sum_{1}^{\infty}|s_{n}|^{-p}<\infty,\end{equation}
 where $s_{n}$ are eigenvalues of the operator $B_{1}$ and $p>1$.
  As operator $B$ is sectorial, then \begin{equation}B_{1}{\geq}0.\end{equation}
\textbf{Theorem 9.3.} \emph{Let the considered Levy process have the
continuous
 density
and let
 the corresponding quasi-potential $B$ be regular and  strongly
 sectorial. If $B_{1}{\in}\Sigma_{p},\quad p>1$ and
\begin{equation}1/p>\beta,\end{equation}
 then the operator $B$ has a point of the
 spectrum different from zero.}\\
\emph{Proof.} It follows from estimation (9.17) that
\begin{equation}||(I-zB)^{-1}||_{p}{\leq}M,\quad
 |\mathrm{arg}z|{\geq}\beta+\epsilon.\end{equation}
  Let us suppose that  the formulated assertion is
not true, i.e. the operator $B$ has no points of the
 spectrum different from zero.
  We set
\begin{equation}
A(r,B)=\mathrm{sup}||(I-re^{i\theta}B)^{-1}||,\quad
0{\leq}\theta{\leq}2\pi.\end{equation} It follows (see [9]) from
condition
 $B_{1}{\in}\Sigma_{p}$
 that $B_{2}{\in}\Sigma_{p}$ and
\begin{equation} \mathrm{log}A(r,B)=O(r^{p}).\end{equation}
According Phragmen -Lindelof
 theorem and relations
 (9.36)-(9.39) we have
\begin{equation}
||(I-zB)^{-1}||{\leq}M.\end{equation}The last relation is possible only when $B=0.$
But in our case $B{\ne}0$. The obtained contradiction proves the theorem.\\
\textbf{Proposition 9.1.} \emph{Let the kernel of $\Phi(x,y)$ of the corresponding
operator
$B$ be bounded. If this operator $B$ is strongly sectorial , then it
has a point of the
 spectrum different from zero.}\\
 \emph{Proof.} As in Theorem 9.3 we suppose that the operator $B$ has no points of the
 spectrum different from zero.
Using the boundedness of the kernel $\Phi(x,y)$
  we obtain the inequality
  \begin{equation} TrB_{1}<\infty.\end{equation} It follows from relations (9.35) and (9.41)
  that (see the triangular model of M. Livshits [15])$\rho=1.$ Since $\beta<1$ all
  the conditions of
  Theorem 9.3 are fulfilled
Hence the proposition is proved.\\
3. Now we shall consider the important case when
\begin{equation}\mathrm{rank}\lambda_{1}=1.\end{equation}
\textbf{Theorem 9.4.} \emph{Let the conditions of theorem $9.2$ be
fulfilled. In the case $(9.42)$ the following relation}
\begin{equation}
 p(t,\Delta)=e^{-t/\lambda_{1}}[c_{1}+o(1)],\quad t{\to}+\infty
 \end{equation}
 \emph{is true.}\\
\emph{Proof.} In view of (9.31)
 we have
\begin{equation}
\mathrm{lim}\frac{1}{T}\int_{0}^{T}q(t)dt{\geq}
|\mathrm{lim}\frac{1}{T}\int_{0}^{T}q(t)e^{-it(\mathrm{Im}{\mu}_{j}^{-1})}dt|,\quad
T{\to}\infty,\end{equation}i.e.
\begin{equation}
g_{1}(0)\int_{\Delta}h_{1}(x)dx{\geq}|\overline{g_{j}(0)}\int_{\Delta}h_{j}(x)dx|.\end{equation}
In the same way we can prove that
\begin{equation}
g_{1}(x_{0})\int_{\Delta_{1}}h_{1}(x)dx{\geq}|\overline{g_{j}(x_{0})}\int_{\Delta_{1}}h_{j}(x)dx|,
\end{equation}where
 \begin{equation}
 x_{0}{\in}\Delta_{1}{\in}\Delta.\end{equation} It follows from (9.46) that
 \begin{equation}
g_{1}(x_{0})h_{1}(x){\geq}|\overline{g_{j}(x_{0})}h_{j}(x)|.
\end{equation}
 We introduce the normalization condition
\begin{equation}g_{1}(x_{0})= g_{j}(x_{0}).\end{equation}
Due to (9.46) and (9.48) the inequalities
\begin{equation}
\int_{\Delta_{1}}h_{1}(x)dx{\geq}|\int_{\Delta_{1}}h_{j}(x)dx|.\end{equation}
 \begin{equation}
h_{1}(x){\geq}|h_{j}(x)|\end{equation} are true. The equality sign
in (9.50) and (9.51) can be only if
\begin{equation}
h_{j}(x)=|h_{j}(x)|e^{i\alpha}.\end{equation} It is possible only in the case when $j=1$. Hence
there exists such a
point $x_{1}$ that
\begin{equation}
h_{1}(x_{1})>|h_{j}(x_{1})|\end{equation} Thus we have
\begin{equation}1=\int_{\Delta_{1}}g_{1}(x)h_{1}(x)dx>
\int_{\Delta_{1}}\overline{g_{j}(x)}h_{j}(x)dx=1,
\end{equation}where $x_{1}{\in}\Delta_{1}.$
 The received contradiction (9.54) means that $j=1.$ Now the assertion
of the theorem follows directly from (9.5) .\\
\textbf{Corollary 9.3.} \emph{Let conditions of Theorem $9.2 $ be
fulfilled. If $\mathrm{rank}\lambda_{1}=1$ and
$x_{0}{\in}\Delta_{1}{\in}\Delta$ then
 the asymptotic equality }\begin{equation}
 p(x_{0},\Delta_{1},t,\Delta)=e^{-t/\lambda_{1}}g_{1}(x_{0})\int_{\Delta_{1}}h_{1}(x)dx[1+o(1)],\quad t{\to}+\infty
 \end{equation}\emph{ is true.}\\The following Krein-Rutman theorem
  [13]  gives the sufficient conditions when relation (9.42) is valid.\\
 \textbf{Theorem 9.5.} \emph{Suppose that the non-negative kernel
 $\Phi(x,y)$ satisfies the condition
 \begin{equation}\int_{\Delta}\int_{\Delta}|\Phi(x,y)^{2}|dxdy<\infty\end{equation}
and has the following property: for each $\epsilon>0$ there exists
 an integer $N=N(\epsilon)$  such that the kernel $\Phi_{N}(x,y)$ of
 a
 the operator $B^{N}$ takes the value zero on a set of points of
 measure not greater than $\epsilon.$ Then}
 \begin{equation}
 \mathrm{rank}\lambda_{1}=1;\quad \lambda_{1}>\lambda_{k},\quad
 k=2,3,...\end{equation}
It is easy to see that the following assertion is valid.\\
 \textbf{Proposition 9.2.} \emph{Let the  inequality
 \begin{equation}\Phi(x,y)>0,\end{equation} be true, when $x{\ne}a_{k},\quad
 x{\ne}b_{k},\quad y{\ne}a_{k},\quad
 y{\ne}b_{k}.$\\
  Then \begin{equation}g_{1}(x)>0 \end{equation} , when} $x{\ne}a_{k},\quad
 x{\ne}b_{k}.$\\
 4.Let us  consider separately the case when the operator $B$ is
 regular and
 \begin{equation}k(x)=k(-x).\end{equation}
 The corresponding operator $S$ is self-adjoint. Hence the operator
 $B$ is self-adjoint and strongly sectorial. In this case equality (9.11) can be
 written in the form
 \begin{equation}
p(t,\Delta)=\sum_{k=1}^{\infty}e^{-t/{\lambda_{k}}}g_{k}(0)\int_{\Delta}g_{k}(x)dx.
\end{equation}
\section{Stable Processes, Main Notions}
1. Let $X_{1},X_{2}...$ be mutually independent random variables
with the same law of distribution $F(x)$. The distribution $F(x)$
is called \emph{strictly stable} if the random variable
\begin{equation}
X=(X_{1}+X_{2}+...+X_{n})/n^{1/\alpha} \end{equation} is also
distributed according to the law $F(x)$. The number $\alpha \quad
(0<\alpha{\leq}2)$ is called a \emph{characteristic exponent} of
the distribution.  The
homogeneous process $X(\tau)\quad (X(0)=0)$ with independent increments
is called a stable process if
\begin{equation}E[\exp{(i{\xi}X(\tau))}]=\exp{\{-\tau|\xi|^{\alpha}[1-i\beta({{\mathrm{sign}}\xi})(\tan{\frac{\pi\alpha}{2}})]\}},
\end{equation} where $0<\alpha<2, \alpha{\ne}1,-1{\leq}\beta{\leq}1,\quad {\tau}>0$.
When $\alpha=1$ we have
\begin{equation}E[\exp{(i{\xi}X(\tau))}]=\exp{\{-\tau|\xi|[1+\frac{2i\beta}{\pi}({{\mathrm{sign}}\xi})(\log{|\xi|)}]\}},
\end{equation} where $-1{\leq}\beta{\leq}1,\quad {\tau}>0$. The
stable processes are a natural generalization of the Wiener
processes. In many theoretical and applied problems it is important
to estimate the value \begin{equation}
p_{\alpha}(t,a)=P(\sup{|X(\tau)|}<a),\quad
0{\leq}\tau{\leq}t.\end{equation} For the stable processes Theorem
9.1. was proved before (see [19]-[22]).
 The value of $p_{\alpha}(t,a)$ decreases very quickly by the
exponential law when $t{\to}\infty$. This fact prompted the idea
to consider the case when the value of $a$ depends on $t$
and $a(t){\to}\infty,\quad t{\to}\infty$. In this paper we deduce
the conditions under which one of the following three cases is
realized:\\
1) $\lim{p_{\alpha}(t,a(t))}=1,\quad, t{\to}\infty.$ \\
2) $\lim{p_{\alpha}(t,a(t))}=0, \quad t{\to}\infty.$\\
3) $\lim{ p_{\alpha}(t,a(t))}=p_{\infty}, \quad  0<p_{\infty}{\leq}1, \quad t{\to}\infty.$\\
We investigate the situation when $t{\to}0$ too.\\
\textbf{Remark 10.1.} In the famous work by M.Kac [11] the
connection of the theory of stable processes and the theory of
integral equations was shown. M.Kac considered in detail only the
case $\alpha=1,\quad \beta=0.$  The case $0<\alpha<2,\quad \beta=0$
was later studied by H.Widom [30]. As to the general case
$0<\alpha<2,\quad -1{\leq}\beta{\leq}1 $ it was investigated in our
works [19]-[22]. In all the mentioned works the parameter $a$ was
fixed. Further we consider the important case when $a$ depends on
$t$ and $a(t){\to}\infty,\quad t{\to}\infty$.
\section{Stable Processes, Quasi-potential.}
1. In this section we  formulate some results from our paper [20]
(see also [22], Ch.7). Here $\psi_{\alpha}(x,s,a)$ is defined by the
relation
\begin{equation}
\psi_{\alpha}(x,s,a)=(I+sB_{\alpha}^{\star})^{-1}\Phi_{\alpha}(0,x,a),\end{equation}
The quasi-potential $B_{\alpha}$ and its kernel $\Phi_{\alpha}(x,y,a)$
will be written later in the explicit form.\\Further we
consider the three cases.\\
\emph{Case 1.} $0<\alpha<2,\quad \alpha{\ne}1,\quad -1<\beta
<1.$\\
\emph{Case 2.} $1<\alpha<2,\quad  \beta={\pm}1.$\\
\emph{Case 3.} $\alpha=1,\quad \beta=0.$\\
Now we introduce the operators \begin{equation} B_{\alpha}f=
\int_{-a}^{a}\Phi_{\alpha}(x,y,a)f(y)dy \end{equation} acting in the
space $L^{2}(-a,a)$. \\In \emph{case 1} the kernel
$\Phi_{\alpha}(x,y,a)$ has the following form
(see[20],[22])\begin{equation}
\Phi_{\alpha}(x,y,a)=C_{\alpha}(2a)^{\mu-1}\int_{a|x-y|}^{a^{2}-xy}{[z^{2}-a^{2}(x-y)^{2}]}^{-\rho}[z-a(x-y)]^{2\rho-\mu}dz,\end{equation}
where the constants $\mu, \rho,$ and $C_{\alpha}$ are defined by the
relations $\mu=2-\alpha,$ \begin{equation}
\sin{\pi\rho}=\frac{1-\beta}{1+\beta}\sin{\pi(\mu-\rho)},\quad
0<\mu-\rho<1,\end{equation}\begin{equation}C_{\alpha}=
\frac{\sin\pi{\rho}}{(\sin{\pi\alpha/2})(1-\beta)\Gamma(1-\rho)\Gamma(1+\rho-\mu)}.\end{equation}
Here $\Gamma(z)$ is Euler's gamma function. We remark that the
constants $\mu, \rho,$ and $C_{\alpha}$ do not depend on
parameter $a$.\\
In \emph{case 2} when $\beta=1$ the relation [20],[22]
\begin{equation} \Phi_{\alpha}(x,y,a)=\frac
{(\cos{\pi\alpha/2})}{(2a)^{\alpha
-1}\Gamma(\alpha)}\{[a(|x-y|+y-x)]^{\alpha-1}-(a-x)^{\alpha-1}(a+y)^{\alpha-1}\}
\end{equation}holds.
In \emph{case 2} when $\beta=-1$ we have [20],[22]
\begin{equation} \Phi_{\alpha}(x,y,a)=\frac
{(\cos{\pi\alpha/2})}{(2a)^{\alpha-1}\Gamma(\alpha)}\{[a(|x-y|+x-y)]^{\alpha-1}-(a+x)^{\alpha-1}(a-y)^{\alpha-1}\}
\end{equation}
Finally, in \emph{case 3} according to M.Kac [11] the equality
\begin{equation}
\Phi_{1}(x,y,a)=\frac{1}{4}{\mathrm{log}}\frac{[a^{2}-xy+\sqrt{(a^{2}-x^{2})(a^{2}-y^{2})}]}{[a^{2}-xy-\sqrt{(a^{2}-x^{2})(a^{2}-y^{2})}]}
\end{equation}is valid.\\
 The  important assertion(see [22],Ch.7) follows from formulas
 (11.2)-(11.8):\\
 \textbf{Proposition 11.1} \emph{Let one of the following conditions be
fulfilled:}\\
I. $0<\alpha<2,\quad \alpha{\ne}1,\quad  -1<\beta<1.$\\
II. $1<\alpha<2,\quad \beta={\pm}1.$\\
III. $\alpha=1, \quad \beta=0.$\\
\emph{Then the corresponding operator $B_{\alpha}$ is regular and strongly sectorial}\\
2. Let us introduce the denotation
\begin{equation}p_{\alpha}(t,,-b,a)=P(-b<X(\tau)<a),\end{equation}
where $a>0,\quad b>0,\quad 0{\leq}\tau{\leq}t.$ We consider in
short the case when  the parameter $b$ is not necessary equal to
$a$. As in case $(-a,a)$ we have the relation
 \begin{equation}
\int_{0}^{\infty}e^{-su}p_{\alpha}(u,,-b,a)du=\int_{-b}^{a}\psi_{\alpha}(x,s,-b,a)dx.\end{equation}
Here $\psi_{\alpha}(x,s,-b,a)$ is defined by relation
\begin{equation}
\psi_{\alpha}(x,s,-b,a)=(I+sB_{\alpha}^{\star})^{-1}\Phi_{\alpha}(0,x,-b,a),\end{equation}
Now the operator $B_{\alpha}$ has the form \begin{equation}
B_{\alpha}f= \int_{-b}^{a}\Phi_{\alpha}(x,y,-b,a)f(y)dy
\end{equation} and acts in the space $L^{2}(-b,a)$. The kernel
$\Phi_{\alpha}(x,y,-b,a)$ is connected with $\Phi_{\alpha}(x,y,a)$
(see (11.7) by the formula
\begin{equation}\Phi_{\alpha}(x,y,-b,a)=
\Phi_{\alpha}(x+\frac{b-a}{2},y+\frac{b-a}{2},\frac{a+b}{2}).\end{equation}
In this way we have reduced the non-symmetric case  $(-b,a)$ to the
symmetric one $(-\frac{a+b}{2},\frac{a+b}{2})$. Let us consider
separately the case $0<\alpha<2,\quad \beta =0.$ In this case the
operator $B_{\alpha}$   is self-adjoint.We denote by
$\lambda_{j},\quad (j=1,2,...)$ the eigenvalues of $B_{\alpha}$ and
by $g_{j}(x)$ the corresponding real normalized eigenfunctions. Then
we can write the new formula (see [20]) for $p_{\alpha}(t,-b,a)$
which is different from (9.10):
\begin{equation}
 p_{\alpha}(t,-b,a)=\sum_{j=1}^{\infty}g_{j}(0)\int_{-b}^{a}g_{j}(x)dxe^{-t\mu_{j}},
 \end{equation}where $\mu_{j}=1/\lambda_{j}$.
\section{On sample functions behavior of stable processes}
From the scaling property of the stable processes we deduce the relations
\begin{equation}p_{\alpha}(t,a)=p_{\alpha}(\frac{t}{a^{\alpha}},1),\end{equation}
\begin{equation}\lambda_{k}(a,\alpha)=a^{\alpha}\lambda_{k}(1,\alpha).\end{equation}
We introduce the notations
\begin {equation}\lambda_{\alpha}(1)=\lambda_{\alpha},\quad
p_{\alpha}(t,1)=p_{\alpha}(t),\quad
g_{\alpha}(x,1)=g_{\alpha}(x),\quad
h_{\alpha}(x,1)=h_{\alpha}(x).\end{equation}
Using relations (12.1), (12.2)
and notations (12.3) we can rewrite Theorem 9.1 in the following
way.\\
\textbf{Theorem 12.1.} \emph{Let one of the following conditions be
fulfilled:}\\
I. $0<\alpha<2,\quad \alpha{\ne}1,\quad  -1<\beta<1.$\\
II. $1<\alpha<2,\quad \beta={\pm}1.$\\
III. $\alpha=1, \quad \beta=0.$\\
\emph{Then the asymptotic equality holds}
\begin{equation}
p_{\alpha}(t,a)=e^{-t/[a^{\alpha}\lambda_{\alpha}]}g_{\alpha}(0)\int_{-1}^{1}h_{\alpha}(x)dx[1+o(1)],
\quad t{\to}\infty.\end{equation}
\emph{Proof.} The corresponding operator $B_{\alpha}$ is regular and strongly
sectorial(see Proposition 11.1).
 The stable processes have the continuous density
(see [31]). So all conditions of Theorem 9.1. are fulfilled. It proves the theorem.\\
\textbf{Remark 12.1.} The operator
$B_{\alpha}$ is self-adjoint when $\beta=0$. In this case
$h_{\alpha}=g_{\alpha}$.\\
\textbf{Remark 12.2.} The value $\lambda_{\alpha}$ characterizes
how fast $p_{\alpha}(t,a)$ converges to zero when $t{\to}\infty$.
The two-sided estimation for $\lambda_{\alpha}$
when $\beta=0$ is given in [17] (see also [22],p.150).\\
3. Now we consider the case when the parameter $a$ depends on $t$.
From Theorem 12.1.  we deduce the assertions.\\
\textbf{Corollary 12.1.} \emph{Let one of conditions I-III of
Theorem $12.1$ be fulfilled and}\begin{equation}
\frac{t}{a^{\alpha}(t)}{\to}\infty,\quad
t{\to}\infty.\end{equation}\emph{Then the following equalities are
true:}\begin{equation} 1)\quad
p_{\alpha}(t,a(t))=e^{-t/[a^{\alpha}(t)\lambda_{\alpha}]}g_{\alpha}(0)\int_{-1}^{1}h_{\alpha}(x)dx[1+o(1)],
\quad t{\to}\infty.\end{equation}
\begin{equation}
2)\quad \lim{p_{\alpha}(t,a)}=0,\quad t{\to}\infty.\end{equation}
\begin{equation}3)\quad \lim{ P[\sup{|X(\tau)|}}>a(t)]=1,\quad
0{\leq}\tau{\leq}t,\quad t{\to}\infty. \end{equation}
\textbf{Corollary 12.2.} \emph{Let one of conditions I-III of
Theorem 12.1 be fulfilled and}\begin{equation}
\frac{t}{[a(t)]^{\alpha}}{\to}0,\quad
t{\to}0.\end{equation}\emph{Then the following equalities are
true:}\begin{equation} 1)\quad \lim{{p_{\alpha}(t,a(t))}}=1,\quad
t{\to}0.\end{equation}
\begin{equation} 2)\quad \lim {P[\sup{|X(\tau)|}}>a(t)]=0\quad
0{\leq}\tau{\leq}t,\quad t{\to}0.
\end{equation}
Corollary 12.2 follows from (12.1) and the relation \begin{equation}
\lim{p_{\alpha}(t)}=1 ,\quad t{\to}0 .\end{equation}
\textbf{Corollary 12.3.} \emph{Let one of conditions I-III of
Theorem 12.1 be fulfilled and}\begin{equation}
\frac{t}{[a(t)]^{\alpha}}{\to}T,\quad 0<T<\infty,\quad
t{\to}\infty.\end{equation}\emph{Then the following equality is
true:}\begin{equation}
\lim{{p_{\alpha}(t,a(t))}}=p_{\alpha}(T),\quad
t{\to}\infty.\end{equation} Corollary 12.3 follows from (12.1).
\section{Wiener Process}
1. We consider separately the important special case when $\alpha
=2$ (Wiener process).  In this case the kernel $\Phi_{2}(x,t,-b,a)$
of the operator $B_{2}$ coincides with the Green's function  (see
[4], [11]) of the equation \begin{equation} -\frac{1}{2}\frac
{d^{2}y}{dx^{2}}=f(x),\quad -b{\leq}x{\leq}a
\end{equation} with the boundary conditions
\begin{equation}
y(-b)=y(a)=0,\quad b>0,\quad a>0. \end{equation}It is easy to see
that
\begin{equation}
\Phi_{2}(x,t,-b,a)=\frac{2}{a+b}\begin{cases}(t+b)(a-x)
&\text{$-b{\leq}t{\leq}x{\leq}a$}\\
(a-t)(b+x) &\text{$-b{\leq}x{\leq}t{\leq}a$}
\end{cases}\end{equation}  Equality (12.1)
is also true when $\alpha=2$ and when $b=a$ , i.e.
\begin{equation}
p_{2}(t,a)=p_{2}(t/a^{2},1). \end{equation} \\
The eigenvalues of problem (13.1),(13.2) have the form
\begin{equation}
\mu_{n}=(\frac{n\pi}{a+b})^{2}/2,\quad n=1,2,3...\end{equation}
The corresponding normalized eigenfunctions are defined by the
equality
\begin{equation} g_{n}(x)=\sqrt\frac{2}{a+b}\sin{[(\frac{n\pi}{a+b})(x+b)]}.\end{equation}
Using formulas (13.5) and (13.6) we have \begin{equation}
p_{2}(t,-b,a)=\sum_{m=0}^{\infty}\frac{4}{(2m+1)\pi}\mathrm{sin}\frac{(2m+1)b\pi}{a+b}e^{-t(\frac{(2m+1)\pi}{a+b})^{2}/2}
\end{equation}
\textbf{Remark 13.1.} If $b=a=1$ then  relation (13.7)takes the form
\begin{equation}
 p_{2}(t)=\sum_{m=0}^{\infty}(-1)^{m}\frac{2}{\pi(m+1/2)}e^{-t[\pi(m+1/2)]^{2}/2}.
 \end{equation} Series (13.8) satisfies the conditions of Leibniz
 theorem. It means that  $p_{2}(t,a)$  can be
 calculated with a given precision when the parameters $t$ and $a$ are
 fixed.\\  From (13.4) and
(13.8) we deduce that \begin{equation} p_{2}(t,a)=
\frac{4}{\pi}e^{-t\pi^{2}/8[a(t)]^{2}}(1+o(1)),
\end{equation} where $t/[a(t)]^{2}{\to}\infty$.\\
\textbf{Proposition 13.1.} \emph{Theorem $12.1$ and Corollaries
$12.1-12.3$
 are true in the
case when }$\alpha=2$ too.\\
\textbf{Remark 13.2.} From the probabilistic point of view  it is
easy to see that the function $p_{2}(t) , (t>0)$ is monotonic
decreasing and \begin{equation} 0<p_{2}(t){\leq}1;\quad\lim
{p_{2}(t)}=1,\quad t{\to}0.\end{equation}
 2. Now we shall describe the behavior of $p(t,-b,a)$ when
 $b{\to}\infty.$ To do it we consider
 \begin{equation}
\frac{d}{dt}p_{2}(t,-b,a)=-\frac{2\pi}{(a+b)^{2}}\sum_{m=0}^{\infty}(2m+1)\mathrm{sin}\frac{(2m+1)b\pi}{a+b}e^{-t(\frac{2m+1}{a+b}\pi)^{2}/2}.
\end{equation}
We use the following Poisson result (see[7]).\\
\textbf{Theorem 13.1.} \emph{If the function $F(x)$ satisfies the
inequalities \begin{equation}
\int_{0}^{\infty}|F(x)|dx<\infty,\quad
\int_{0}^{\infty}|F^{\prime}(x)|dx<\infty \end{equation} then the
equality
\begin{equation}
\sum_{m=0}^{\infty}F(m)=\frac{1}{2}F(0)+\int_{0}^{\infty}F(x)dx+
2\sum_{m=1}^{\infty}\int_{0}^{\infty}F(x)\mathrm{cos}2{\pi}mxdx
\end{equation}is true.}\\
Thus in case (13.11) we have \begin{equation}
F(x)=G(x)-G(2x),\end{equation} where \begin{equation}G(x)=
-\frac{2\pi}{(a+b)^{2}}x\mathrm{sin}\frac{xb\pi}{a+b}e^{-t(\frac{x}{a+b}\pi)^{2}/2}.
\end{equation}It is easy to see that conditions (13.12) are
fulfilled and \begin{equation} F(0)=0,\quad
\int_{0}^{\infty}F(x)dx=\frac{1}{2}\int_{0}^{\infty}G(x)dx
\end{equation}Using (13.15) and (13.16) we deduce the equality
\begin{equation}\int_{0}^{\infty}F(x)dx=-\frac{1}{{\pi}t}\int_{0}^{\infty}ue^{-u^{2}/2}\mathrm{sin}\frac{ua}{\sqrt{t}}du,
\end{equation}where $u=\frac{x\pi}{a+b}\sqrt{t}$. Now we use the
following relation from the sine transformation theory
(see[25]).\begin{equation}
\int_{0}^{\infty}ue^{-u^{2}/2}\mathrm{sin}(xu)du=\sqrt{\frac{\pi}{2}}xe^{-u^{2}/2}.\end{equation}
In view of (13.17) and  (13.18) the equality
\begin{equation}
\int_{0}^{\infty}F(x)dx=-\frac{a}{\sqrt{2\pi}}t^{-3/2}e^{-a^{2}/2t}\end{equation}
is true. Now  we calculate the integrals
\begin{equation}
J_{m}=2\int_{0}^{\infty}G(2x)\mathrm{cos}2{\pi}mxdx,\quad
I_{m}=2\int_{0}^{\infty}G(x)\mathrm{cos}2{\pi}mxdx
\end{equation}
Using again formula (13.18) we have
\begin{equation}
J_{m}=-\sqrt{2/\pi}t^{-3/2}[A_{m}e^{-A_{m}^{2}/2t}-B_{m}e^{-B_{m}^{2}/2t}],
\end{equation}where $A_{m}=2m(a+b)+a,\quad B_{m}=2m(a+b)-a.$ In
the same way we found \begin{equation}
I_{m}=-\sqrt{1/2\pi}t^{-3/2}[C_{m}e^{-C_{m}^{2}/2t}-D_{m}e^{-D_{m}^{2}/2t}],
\end{equation}where $C_{m}=m(a+b)+a,\quad D_{m}=m(a+b)-a.$ From
relation (13.7) and equality
\begin{equation}
[\int_{c/\sqrt{t}}^{d/\sqrt{t}}e^{-u^{2}/2}du]^{\prime}=-\frac{1}{2}t^{-3/2}(de^{-d^{2}/2t}-ce^{-c^{2}/2t})
\end{equation}we obtain the following representation of
$p_{2}(t,-b,a)$:
\begin{equation}
p_{2}(t,-b,a)=1-\sqrt{2/{\pi}}\int_{a/\sqrt{t}}^{\infty}e^{-u^{2}/2t}du+q_{\alpha}(t,-b,a),\end{equation}
where \begin{equation}
q_{2}(t,-b,a)=\sqrt{2/{\pi}}\sum_{m=1}^{\infty}[2\int_{B_{m}/\sqrt{t}}^{A_{m}/\sqrt{t}}e^{-u^{2}/2t}du-
\int_{D_{m}/\sqrt{t}}^{C_{m}/\sqrt{t}}e^{-u^{2}/2t}du]\end{equation}
So we have deduced two formulas (13.7) and (13.24) for
$p_{2}(t,-b,a)$. Formula (13.7) is useful when $t$ is great and
the parameters $a$ and $b$ are fixed. \\
\textbf{Proposition 13.2.} \emph{In the case of the Wiener process ($\alpha=2$)
the asymptotic equality}
\begin{equation}
p_{2}(t,-b,a)=\frac{4}{\pi}\mathrm{sin}\frac{a\pi}{a+b}e^{-t(\pi)^{2}/[2(a+b)^{2}]}[1+o(1)],
\quad t{\to}\infty
\end{equation} \emph{holds}\\
Formula (13.24) is useful when $b$ is great and parameters
$a$ and $t$ are fixed. \\
\textbf{Proposition 13.3.} \emph{In the case of the Wiener process ($\alpha=2$)
the asymptotic equality}
$p_{2}(t,-b,a)=$
\begin{equation}
=1-\sqrt{2/{\pi}}\int_{a/\sqrt{t}}^{\infty}e^{-u^{2}/2t}du-\sqrt{2/{\pi}}\int_{b/\sqrt{t}}^{(b+2a)/\sqrt{t}}e^{-u^{2}/2t}du[1+o(1)],
\end{equation}\emph{where $b{\to}\infty$, is valid.}\\
The well-known formula (see [8]) for \emph{the first hitting time}
\begin{equation}
p_{2}(t,-\infty,a)=1-\sqrt{2/{\pi}}\int_{a/\sqrt{t}}^{\infty}e^{-u^{2}/2t}du
\end{equation} follows directly from (13.27).
\section{Iterated logarithm law, most visited sites and  first
hitting time} It is interesting to compare our results (Theorem
9.1,
 Corollaries 12.1-12.3 and Proposition 13.1 -13.3) with the well-known results
  mentioned in the title of
 the section.\\
1. We begin with the famous Khinchine   theorem (see[4]) about the
iterated
logarithm law.\\
\textbf{Theorem 14.1.} \emph{Let $X(t)$ be stable process
($0<\alpha<2$). Then almost surely (a.s.) that} \begin{equation}
\lim
{\frac{\sup{|X(t)|}}{(t\log{t})^{1/\alpha}|\log{|\log{t}|}|^{(1/{\alpha)+\epsilon}}}}=
\begin{cases}0 &\text{$\epsilon>0$ a.s.}\\
\infty &\text{$\epsilon=0$ a.s.}
\end{cases}\end{equation}
 We introduce the random process \begin{equation}
U(t)=\sup{|X(\tau)|}, \quad {0\leq}\tau{\leq}t \end{equation}
From Corollaries 12.1-12.3 and Proposition 13.1 we deduce the assertion.\\
\textbf{Theorem 14.2.} \emph{Let one of  conditions $I-III$ of
Theorem $12.1$ be fulfilled or let $\alpha=2$ and
\begin{equation} b(t){\to}\infty ,\quad t{\to}\infty.\end{equation} Then }
\begin{equation} b(t)U(t)/t^{1/\alpha}{\to}\infty\quad (P),\quad U(t)/[b(t)t^{1/\alpha}]{\to}0 \quad (P)\end{equation}\\
(It is denoted by symbol $(P)$, that the
convergence is in probability .)\\
In particular we have :
\begin{equation} [(\log^{\epsilon}{t})U(t)]/t^{1/\alpha}{\to}\infty\quad (P),
\quad U(t)/[(\log^{\epsilon}{t})t^{1/\alpha}]{\to}0 \quad (P),\end{equation}\\
when $\epsilon>0.$ We see that our approach and the
classical one have some similar points (estimation of
$|X(\tau)|$), but these approaches are essentially different. We
consider the behavior of $|X(\tau)|$ on the interval $(0,t)$, and
in the classical case $|X(\tau)|$
is considered on the interval $(t,\infty)$.\\
2.  We denote by $V(t)$ the most visited site of stable process $X$
up to time $t$ (see [1]). We formulate the following result
(see [1] and references therein).\\
Let $1<\alpha<2,\quad \beta =0,\quad \gamma
>9/(\alpha-1)$. Then the relation \begin{equation}
\lim{\frac{(\log{t})^{\gamma}}{t^{1/\alpha}}|V(t)|}=\infty ,\quad
t{\to}\infty\quad (a.s.)\end{equation}is true.
 To this important result we add the following estimation.\\
\textbf{Theorem 14.3.} \emph{Let one of the conditions $I-III$ of
Theorem $12.1$ be fulfilled or let $\alpha=2$ and
\begin{equation} b(t){\to}\infty ,\quad t{\to}\infty.\end{equation} Then }
\begin{equation} |V(t)(t)/[b(t)t^{1/\alpha}]{\to}0 \quad (P)\end{equation}\\
In particular we have when $\epsilon>0$:
\begin{equation}  |V(t)|/[(\log^{\epsilon}{t})t^{1/\alpha}]{\to}0 \quad (P)\end{equation}\\
 The formulated theorem follows directly from the
inequality $U(t){\geq}|V(t)|$. \\
3. The first hitting time $T_{a}$ is defined by the formula
\begin{equation}T_{a}=\inf{(t{\geq}0,X(t){\geq}a)}.
\end{equation} It is obvious that
\begin{equation}P(T_{a}>t)=P[\sup{X(\tau)}<a],\quad
0{\leq}\tau{\leq}t.\end{equation} We have
\begin{equation}P(T_{a}>t){\geq}P[-b<\sup{X(\tau)}<a]=p_{\alpha}(t,-b,a),\quad
0{\leq}\tau{\leq}t.\end{equation}So our formulas for
$p(t,-b,a)$ estimate $P(T_{a}>t)$ from below. It is easy
to see that \begin
{equation}p(t,-b,a){\to}P(T_{a}>t),\quad
b{\to}+\infty.\end{equation}
\\\textbf{Remark 14.1.} Our results can be interpreted in terms of
the  first hitting time $T_{[-b,a]}$ one of the barriers either $-b$
or $a$ (ruin problem). Namely, we have
\begin{equation} P(T_{[-b,a]}>t)=p(t,-b,a).\end{equation}
The distribution of the first hitting time for the Levy
processes
 is an open problem.\\
 \textbf{Remark
14.2.} Rogozin B.A. in his interesting work [18] established the law
of the overshoot distribution for the stable processes when the
existing interval is fixed.\\

\begin{center}\textbf{References}\end{center}
1. \textbf{Bass R.F., Eisenbaum N. and Shi Z.}, The Most Visited
Sites of Symmetric Stable Processes,  Probability Theory and
Related Fields, 116,
(2000), 391-404.\\
2.\textbf{Bonsall F.F., Duncan J.,} Numerical Ranges.1-49, MAA
Studies in Mathematics, v.21 (ed. Bartle R.G.) 1980.\\
3.\textbf{Baxter G., Donsker M.D.,} On the Distribution of the
Supremum Functional for Processes with Stationary Independent
Increments, Trans. Amer. Math. Soc. 8,
73-87,1957.\\
4. \textbf{Bertoin J.,} Levy Processes, University Press, Cambridge, 1996.\\
5. \textbf{Chuangyi Z.,} Almost Periodic Type Functions and
Ergodicity,
Beijing, New York, Kluwer, 2003.\\
6. \textbf{Chung K.L.,} Green, Brown and Probability , World
Scietific, 2002.\\
7. \textbf{Evgrafov M.A.} Asymptotic Estimates and Entire Functions
, Gordon and Breach,
New York, 1961.\\
8. \textbf{Feller W.}, An Introduction to Probability Theory and its
Applications, J.Wiley and Sons, 1971.\\
9.\textbf{Gohberg I., Krein M.G.,} Introduction to the Theory of
Non-selfadjont Operators , Amer. Math. Soc. Providence, 1970.\\
10. \textbf{Ito K.}, On Stokhastic Differential Equations, Memoirs
Amer. Math.
 Soc. No.4, 1951.\\
11. \textbf{Kac M.}, On some Connections Between Probability Theory
and Differential and Integral Equations, Proc.Sec.Berkeley
Symp.Math.Stat. and Prob., Berkeley, 189-215, 1951.\\
12. \textbf{Kac M.}, Probability and Related Topics in Physical Sciences, Colorado, 1957.\\
13. \textbf{Krein M.G.,Rutman M.A.,} Linear Operators Leaving
Invariant a Cone
in a Banach Space, Amer. Math. Soc. Translation, no.26, 1950.\\
14. \textbf{Levitan B.M.} Some Questions of the Theory  of Almost
Periodic Functions,
Amer. Math. Soc.,Translation, 28, 1950.\\
15. \textbf{Livshits M.S.,} Operators, Oscillations, Waves, Open
Systems, American Math.
 Society, Providence, 1973.\\
16.\textbf{Pietsch A.,}  Eigenvalues and s-Numbers, Cambridge
University Press, 1987.\\
17. \textbf{Pozin S.M., Sakhnovich L.A.,} Two-sided Estimation of
the Smallest Eigenvalue of an Operator Characterizing Stable
Processes, Theory Prob. Appl., 36, No.2,385-388,
1991 \\
18. \textbf{Rogozin B.A.,}The distribution of the first hit for
stable and asymptotically stable walks on an interval, Theory
Probab. Appl. 17, 332-338, 1972.\\
19. \textbf{Sakhnovich L.A.,} Abel Integral Equations in the theory
of Stable Processes, Ukr.Math. Journ., 36:2, 193-197 ,
1984.\\
20. \textbf{Sakhnovich L.A.,} Integral Equations  in the theory of
Stable Processes, St.Peterburg Math.J.,4, No.4 1993,819-829.\\
21. \textbf{Sakhnovich L.A.,} The Principle of Imperceptibility of
the Boundary in the Theory of Stable Processes , St.Peterburg
Math.J.,
6, No.6, 1995, 1219-1228.\\
22. \textbf{Sakhnovich L.A.,} Integral Equations with Difference
Kernels, Operator Theory, v.84, 1996, Birkhauser.\\
23. \textbf{Sato K.}, Levy Processes and Infinitely Divisible
Distributions, University
Press, Cambridge, 1999.\\
24. \textbf{Schoutens W.,} Levy Processes in Finance, Wiley series
in
Probability and Statistics,2003\\
25. \textbf{Titchmarsh E.C.}, Introduction to the Theory of Fourier Integrals,Oxford, 1937.\\
26.\textbf{Stone M.,} Linear Transformation in Hilbert space., New York, 1932.\\
27. \textbf{Thomas M., Barndorff O.,}(ed.), Levy Processes ;
Theory and Applications, Birkhauser, 2001.\\
28. \textbf{Trefethen L.N., Embree M.} Spectra and Pseudospectra ,
 Princeton University
Press,2005.\\
29. \textbf{Tuominen P., Tweedie R.L.}, Exponential Decay and
Ergodicity  of General Markov
Processes and their Discrete Skeletons. Adv.in Appl. Probab. 11, 784-803.1979.\\
30. \textbf{Widom H.}, Stable Processes and Integral Equations,
Trans. Amer. Math. Soc. 98, 430-449, 1961.\\
31. \textbf{Zolotarev V.M.}, One-dimensional stable distribution,
Providence, Amer. Math. Soc. 1986.\\
\end{document}